\documentclass[11pt,reqno]{amsart}

\usepackage[numbers,sort&compress]{natbib}
\usepackage[colorlinks,citecolor=blue]{hyperref}
\usepackage{amssymb}

\usepackage{enumerate}

\usepackage{graphicx}
\usepackage{mathrsfs}
\usepackage{graphicx}
\usepackage{subfigure}
\usepackage{epstopdf}
\usepackage{amssymb}
\usepackage{caption}
\usepackage{subfigure}
\usepackage{verbatim}
\usepackage{color}
\usepackage{amsmath}
\usepackage{cases}
\usepackage{float}
\makeatletter\@addtoreset{equation}{section}

\makeatletter
\@namedef{subjclassname@2010}{%
	\textup{2010} Mathematics Subject Classification}
\makeatother

\newtheorem{thm}{Theorem}[section]

\theoremstyle{definition}

\newtheorem{rem}[thm]{Remark}

\numberwithin{equation}{section}

\makeatletter\@addtoreset{equation}{section}

\newenvironment{proofed}[1]{\par \textbf{Proof}\quad #1}{\hfill \textbf{} $\Box$ }

\newtheorem{example}{\bf{Example}}[section]

\numberwithin{equation}{section}
\numberwithin{table}{section}
\numberwithin{figure}{section}
\newcommand{\ba}{\begin{array}}\newcommand{\ea}{\end{array}}
\newcommand{\be}{\begin{eqnarray}}\newcommand{\ee}{\end{eqnarray}}
\newcommand{\beq}{\begin{equation*}}\newcommand{\eeq}{\end{equation*}}
\newcommand{\bex}{\begin{eqnarray*}}
	\newcommand{\eex}{\end{eqnarray*}}
\newcommand{\bse}{\begin{subequations}}
	\newcommand{\ese}{\end{subequations}}
\newcommand{\bal}{\begin{align}}
\newcommand{\eal}{\end{align}}

\usepackage{lineno}
\usepackage{enumerate}
\font\tenbi=cmmib10   at 11 pt \font\sevenbi=cmmib10 at 9pt
\font\fivebi=cmmib7 at 6pt
\newfam\bifam
\textfont\bifam=\tenbi 
\scriptfont\bifam=\sevenbi
\scriptscriptfont\bifam=\fivebi

\def\bi{\fam\bifam\tenbi}
\font\tendb=msbm10 at 12 pt \font\sevendb=msbm7
\newfam\dbfam
\textfont\dbfam=\tendb \scriptfont\dbfam=\sevendb

\def\n{{\bi n}}
\def\x{{\bi x}}

\def\dps{\displaystyle}

\begin{document}
	\graphicspath{{figure/};}
	\baselineskip=17pt
	
	\title
	[New schemes for anisotropic phase-field crystal growth model]
	{New efficient time-stepping schemes for the anisotropic phase-field dendritic crystal growth model$^\ast$}
	
	\author[M.H. Li]{Minghui Li$^{1}$}
	\author[M. Azaiez]{Mejdi Azaiez$^{1,2}$}
	
	\author[C.J. Xu]{Chuanju Xu$^{1,3}$}
	\thanks{\hskip -12pt
		$^\ast$This research is partially supported by NSFC grant 11971408, NNW2018-ZT4A06 project, and NSFC/ANR joint program 51661135011/ANR-16-CE40-0026-01.\\
		$^{1}$School of Mathematical Sciences and
		Fujian Provincial Key Laboratory of Mathematical Modeling and High Performance
		Scientific Computing, Xiamen
		University, 361005 Xiamen, China.
		\\$^{2}$Institut Polytechnique de Bordeaux, Laboratoire I2M CNRS UMR5295, France.\\
		${}^{3}$Corresponding author. Email: cjxu@xmu.edu.cn (C. Xu)}
	\date{}

	\begin{abstract}
		In this paper, we propose and analyze a first-order and a second-order time-stepping 
		schemes for the anisotropic phase-field dendritic crystal growth model.
		The proposed schemes are based on an auxiliary variable approach 
		for the Allen-Cahn equation and
		delicate treatment of the terms coupling the Allen-Cahn equation and temperature equation. 
		The idea of the former is to introduce suitable auxiliary variables to 
		facilitate construction of high order stable schemes for a large class of gradient flows.
		We propose a new technique to treat the coupling terms involved in the crystal growth model, 
		and introduce suitable stabilization terms to result in  
		totally decoupled schemes, which satisfy a discrete energy law 
		without affecting the convergence order. 
		A delicate implementation demonstrates that the proposed schemes can be realized in a very efficient
		way. That is, it only requires solving four linear elliptic equations and a simple algebraic equation at
		each time step. 
		A detailed comparison with existing schemes is given, and the advantage of the new schemes 
		are emphasized. 
		As far as we know this is the first second-order scheme that is totally decoupled, linear, unconditionally stable for the dendritic crystal growth model with variable mobility parameter.
	\end{abstract}
	
	
	\maketitle

	\section{Introduction}
	On one side, dendritic growth is a very common phenomenon in nature. 
	We are all familiar with the way how trees grow by spreading branches and roots 
	from the main trunk. This is where the name ``dendritic" comes from, although
	the term ``dendrite" itself is used to describe branched projections of neurons.
	On the other side, dendritic growth phenomena and the shapes of growing crystals 
	are of fundamental interest to physicists and are of practical importance to engineers. 
	Crystal dendritic growth is one of the most extensively studied topics 
	in the scientific literature. 
	Crystallization proceeds through the competition between thermodynamics
	-- driven by the local undercooling of the liquid ahead of the solidification front -- and
	the ability of the system to diffuse latent heat of fusion away from the solid-liquid interface.
	It usually forms natural fractal microstructure, so-called dendrites, 
	which are the ubiquitous crystal form in freezing alloys and supercooled melts.
	When the molten material is supercooled below the freezing point of the solid, 
	a spherical solid nucleus grows in the undercooled melt initially. 
	Along with some preferred directions of growth, 
	the solid form begins to express some protrusion accompanied by steeper concentration gradients at its end.
	Dendritic microstructures formed during solidification/freezing play 
	a key role in properties of the final solid material. 
	Understanding these microstructures is therefore considered essential for controlling basic solidification and crystal growth processes. 
	
	The first phase-field models were suggested for numerically simulating dendritic growth
	in 1980s; see, e.g., 
	\cite{Fix83,CL85,Cag86,Lan86}.
	This concept has been validated by comparison
	with theoretical predictions and experimental measurements and is applied to a broad range of 
	investigations in materials science \cite{BWBK02}.
	Nowadays, the phase-field method has emerged as a powerful tool for modelling and simulation of crystal dendritic growth. 
	In contrast to sharp interface approaches with interfaces of zero thickness, 
	the phase-field model introduces a smooth phase-field variable 
	by a diffuse interface profile to distinguish between the solid and liquid phases.
	In this model the complicated topological changes of a solid-liquid interface
	can be handled in an easy way without the need of the explicit tracking of the interface. 
	The phase field is considered as an order parameter which is introduced to describe the moving interfacial
	boundary between unstable and stable phases during phase transformation processes. 
	By asymptotic expansions, it can be shown that the phase-field methods relate to classical sharp interface models 
	such as Hele-Shaw type models and Stefan problems in the limit of zero interfacial thickness, see, e.g. \cite{CS94}.
	
	Another advantage of the phase-field approach is that the governing set of equations in the model
	can be naturally derived from an energy-based variational principle. The variational framework of phase-field formulations makes them thermodynamically consistent and physically attractive in modeling the general phase-field dendritic crystal growth model
	\cite{KR98,WB95,RB05examination,KPGD99universal,MMM96phase,DDFR13dendritic,M03study,DZP17phase}.
	
	In this paper, we will focus on the numerical approximations
	for the anisotropic phase-field dendritic crystal growth model proposed in \cite{KR98,KR99}. The model is composed of two coupled nonlinear equations. One is the phase-field equation that governs the anisotropy of the crystal. 
	The other is the heat equation that controls heat diffusion of the system. 
	It is shown that this nonlinear
	coupled system satisfies a thermodynamically consistent energy dissipation 
	law.
	The main aim of this paper is to design efficient numerical schemes for this
	nonlinear crystal growth model, 
	which satisfies a discrete version of the energy dissipation law. 
	In fact, constructing schemes that preserve the discrete energy dissipation law
	for similar models has been subject of many recent papers
	\cite{Beckerrohl2008,shen2010numerical,fjordholm2011well,wang2011energy,chen2012linear,zhang2013adaptive,chen2013time,shen2015decoupled,cheng2018multiple}. 
	Although large amounts of works have been devoted to numerical approximation
	for phase-field dendritic crystal growth models; see, e.g.,  \cite{ZWY17numerical,LZW2019energy,SHS14,LK12phase} and the references therein,  
	there is still a need for efforts on developing low-cost, stable, and high order schemes for such models. 
	
	The main difficulties in constructing highly efficient schemes for the dendritic crystal growth models come from: 
	1) the double-well energy potential and 
	the stiffness associated with the interfacial width in the phase equation;
	2) the anisotropic coefficient; 
	3) the nonlinear interaction terms in both the heat equation and 
	the phase field equation.
	Let's briefly review recent progress in this direction. Firstly,
	to overcome the difficulty caused by the nonlinearity and the thin interface in the phase field equation,
	schemes based on the invariant energy quadratization (IEQ) \cite{Y2016linear} and scalar auxiliary variable approach (SAV) \cite{shen2018scalar,SY2020ieq} have been proposed: 
	a decoupled stable but only first-order scheme in \cite{ZCY19novel} 
	and a second-order stable scheme but fully coupled scheme in \cite{Y19efficient}.  Due to the presence of the nonlinear phase term in the heat equation and the interaction term 
	in the phase equation, it seems not easy to design a fully-decoupled, second-order accurate and energy stable scheme. 
	For example, 
	the splitting method used in \cite{ZCY19novel} is not directly extendable 
	to a second-order discretization for the time derivative of the phase function 
	in the heat equation.
	In the case of constant mobility parameter,  
	\cite{Y21anovel} proposes a decoupling, linear, second-order accurate, and unconditionally stable scheme by using multi-auxiliary variables.
	A similar technique was used in \cite{Y21new,Y21novel,Y21numerical}
	to deal with some coupling models. 
	However, 
	compared with the traditional SAV approach, multi-auxiliary variable approach means extra computational cost since more equations are to be solved.
	In particular, the second-order scheme proposed in \cite{Y21anovel} is a four-step 
	scheme, thus is much more computationally expensive.
	
	The main purpose of the present paper is to propose easy-to-implement, 
	second-order accurate, and unconditionally stable schemes 
	for the anisotropic phase-field dendritic crystal growth model. First of all, we rewrite the time derivative term of the phase function  
	in the heat equation into an equivalent form, which allows to design a 
	three-step second-order scheme. 
	The idea is to introduce a suitable auxiliary variable 
	to the Allen-Cahn equation and a new technique to treat the coupling terms. 
	Then some carefully chosen stabilization terms are added to result in  
	totally decoupled schemes that satisfy a discrete energy law 
	without losing the convergence order. 
	A careful examination shows that the proposed schemes can be implemented 
	by only solving four linear elliptic equations and a simple algebraic equation. 
	As far as we know this is the first second-order scheme that is totally decoupled, linear, unconditionally stable for the dendritic crystal growth model
	variable mobility parameter.
	In the case of constant mobility parameter, 
	compared with \cite{Y21anovel} (a four-layer scheme that requires solve five linear elliptic equations with constant coefficients and some algebraic equations), our scheme is a three-layer scheme that only needs solve four linear elliptic equations with constant coefficients and a simple algebraic equation.
	
	The rest of the paper is organized as follows. 
	In Section 2, we describe the phase-field dendritic crystal growth model, and 
	present the equivalent reformulation using auxiliary variables.
	In Section 3 we propose a first-order unconditionally stable time-stepping scheme, 
	and prove the energy decay property of the proposed scheme.
	Section 4 is devoted to construct and analyze a second-order, linear, decoupled, and unconditionally stable scheme.
	The implementation detail is also presented to show 
	that the scheme can be efficiently realized through solving a set of decoupled, linear elliptic equations. 
	We give in Section 5 some numerical examples to verify the efficiency of the proposed methods. 
	Finally, the paper ends with some concluding remarks.
	
	\section{Governing equations and auxiliary variable reformulation}
	\setcounter{equation}{0}
	\subsection{Governing equations}
	We are interested in numerically solving the following anisotropic phase-field dendritic crystal growth model equations 
	in the domain $\Omega\subset R^{2}$:
	\be
	&&\varrho(\phi) \phi_{t}=-\frac{\delta E}{\delta \phi}-\frac{\lambda}{\varepsilon} h^{\prime}(\phi) T,\label{problem-3}\\
	&&T_{t}= \nabla \cdot (D \nabla T)+K h^{\prime}(\phi) \phi_{t}, \label{problem-4}
	\ee
	where $\phi(\x, t)$ is the phase function to label the liquid and solid phases, 
	$\varrho(\phi)>0$ is the mobility parameter that can be chosen either as a constant \cite{Y21anovel}, 
	or as a function of $\phi$ \cite{Y19efficient}. 
	$\varepsilon>0$ is a parameter used to control the interface width, 
	$\lambda$ is the linear kinetic coefficient.
	In Equation \eqref{problem-4}, 
	$T(\boldsymbol{x}, t)$ is the scaled temperature, $D$ is the constant diffusion rate of the temperature, and
	$K$ is the latent heat parameter that controls the speed of heat transfer along with the interface. 
	It is worth noting that the efficiency of the schemes we propose below covers the case 
	$D$ is a function of $\phi$ \cite{S2002phase}. 
	The function ${h(\phi)}$ is defined by 
	\bex
	h(\phi) := \frac{1}{5} \phi^{5}-\frac{2}{3} \phi^{3}+\phi,
	\eex 
	which represents a generation of latent heat.   
	Following the phenomenological free energy used in \cite{KR98}, we consider here 
	\begin{equation}\label{E1}
	E(\phi, T)=\int_{\Omega}\left(\frac{1}{2} \kappa^{2}(\nabla \phi)|\nabla \phi|^{2}+\frac{1}{\varepsilon^{2}} F(\phi)+\frac{\lambda}{2 \varepsilon K} T^{2}\right) d \x,
	\end{equation}
	where $F(\phi)=\frac{1}{4}\left(\phi^{2}-1\right)^{2}$ is the double-well type Ginzburg-Landau potential.
	$\kappa(\cdot)$ in \eqref{E1} is a function describing the anisotropic property, 
	which takes the form \cite{KR98,KR99}: 
	\be\label{Ka}
	\kappa(\nabla \phi)=1+ \sigma \cos (m \theta),
	\ee
	where $m$ is a model number  of anisotropy, $\sigma$ is the parameter for the anisotropy strength, and 
	$\theta=\arctan \big(\frac{\phi_{y}}{\phi_{x}}\big)$.
	The variational derivative of $E$ with respect to $\phi$ is:
	$$
	\frac{\delta E}{\delta \phi}=-\nabla \cdot\left(\kappa^{2}(\nabla \phi) \nabla \phi+\kappa(\nabla \phi)|\nabla \phi|^{2} \boldsymbol{H}(\phi)\right)+\frac{f(\phi)}{\varepsilon^{2}},
	$$
	where $\boldsymbol{H}(\phi)$ is the variational derivative of $\kappa(\nabla \phi)$, and $f(\phi)=F^{\prime}(\phi)$.
	In the case $m=4$, a direct calculation shows 
	\be\label{H}
	\boldsymbol{H}(\phi):=\frac{\delta \kappa(\nabla \phi)}{\delta \phi}
	=4 \sigma \frac{4}{|\nabla \phi|^{6}}\Big(\phi_{x}\big(\phi_{x}^{2} \phi_{y}^{2}-\phi_{y}^{4}\big), 
	\phi_{y}\big(\phi_{x}^{2} \phi_{y}^{2}-\phi_{x}^{4}\big)\Big).
	\ee

	For convenience, we only consider that the equations \eqref{problem-3} and \eqref{problem-4} are subject to the 
	Neumann boundary conditions
	\begin{equation}\label{bc}
	\frac{\partial \phi}{\partial \n}|_{\partial \Omega}=0,\quad \frac{\partial T}{\partial \n}|_{\partial \Omega}=0,
	\end{equation} 
	although other boundary conditions such as the periodic conditions are possible.
	
	Now we briefly recall main property of \eqref{problem-3}-\eqref{bc}.
	A key property of the model is that it satisfies an energy law, which can be
	derived by 
	taking the inner product of \eqref{problem-3} with $-\phi_{t}$  and using integration by parts:
	\bex
	\frac{d}{dt}\int_{\Omega}\left(\frac{1}{2} \kappa^{2}(\nabla \phi)|\nabla \phi|^{2}+\frac{1}{\varepsilon^{2}} F(\phi) \right)d \x +\int_{\Omega}\frac{\lambda}{\varepsilon} h^{\prime}(\phi) T \phi_{t} d \x=-\int_{\Omega}\varrho(\phi) \phi^2_{t} d \x.
	\eex
	Then
	taking the inner product of \eqref{problem-4} by $-\frac{\lambda}{\varepsilon K} T$
	gives:
	\bex
	\frac{d}{dt}\int_{\Omega}\frac{\lambda}{2 \varepsilon K} T^{2} d \x-\int_{\Omega}\frac{\lambda}{\varepsilon} h^{\prime}(\phi) T \phi_{t} d \x=-\frac{\lambda D}{\varepsilon K}\int_{\Omega} \nabla T\cdot \nabla T  d \x.
	\eex
	Combining the above two equalities gives the following energy law
	\be\label{energylaw}
	\frac{d}{d t} E(\phi, T)=-\left\|\sqrt{\varrho(\phi)} \phi_{t}\right\|^{2}-\frac{\lambda D}{\varepsilon K}\|\nabla T\|^{2},
	\ee
	where $\|\cdot\|$ denotes the standard $L^{2}(\Omega)$ norm. This means that 
	the energy $E(\phi, T)$ decays in time during the crystal-growing process.
	
	\subsection{Auxiliary variable reformulation}
	
	The main purpose of this paper is to develop novel efficient schemes for the anisotropic crystal growth model \eqref{problem-3}-\eqref{bc}.
	We start with an auxiliary variable approach, 
	which will be used later to construct time-stepping schemes for the phase field 
	equation \eqref{problem-3}. 
	We define the variable 
	\begin{equation}\label{E1B}
	R(t)=\sqrt{E_{1}(\phi)}, \ \ \ E_{1}(\phi)={\int_{\Omega}\left(\frac{1}{2}\big(\kappa^{2}(\nabla \phi)-S_{1}\big)|\nabla \phi|^{2}+\frac{1}{\varepsilon^{2}}\Big(F(\phi)-\frac{S_{2}}{2} \phi^{2}\Big)+B \right)d\x},
	\end{equation}
	where $S_{1}$ and $S_{2}$ are two positive constants, $0<S_{1}<\left(1-\sigma\right)^{2}$,
	$B$ is a positive constant used to make $E_{1}$ positive. 
	Notice $\kappa^{2}(\nabla \phi) \geq\left(1-\sigma\right)^{2}$ and $F(\phi)$ is a quartic polynomial, one can verify that 
	${\int_{\Omega}(\frac{1}{2}\big(\kappa^{2}(\nabla \phi)-S_{1}\big)|\nabla \phi|^{2}+\frac{1}{\varepsilon^{2}}(F(\phi)-\frac{S_{2}}{2} \phi^{2}) )d\x}$ is bounded from below. Therefore such a constant $B$ exists.
	The introduction of the constants $S_1$ and $S_2$ is inspired by the work \cite{Y19efficient,ZY20fully,Y21anovel}.
	We will see that these constants help in ensuring the $H^1$-stability of the phase function. 
	
	Using the auxiliary variable $R(t)$, the total free energy \eqref{E1} can be rewritten as
	\begin{equation}\label{E3}
	E(\phi, R, T)=\int_{\Omega}\left(\frac{\lambda}{2 \varepsilon K} T^{2}+\frac{S_{1}}{2}|\nabla \phi|^{2}+\frac{S_{2}}{2 \varepsilon^{2}} \phi^{2}-B\right) d \boldsymbol{x}+R^2,
	\end{equation}
	and the original equations \eqref{problem-3}-\eqref{problem-4} can be reformulated into the following equivalent form:
	\bse\label{sche11}
	\bal
	& \phi_{t}=M(\phi)\mu,\label{problem-31}\\
	&\mu=-\frac{R(t)}{\sqrt{E_{1}(\phi)}} g(\phi)+S_1\Delta\phi-\frac{S_2}{\varepsilon^2}\phi-\frac{R(t)}{\sqrt{E_{1}(\phi)}}\frac{\lambda}{\varepsilon} h^{\prime}(\phi) T,\label{problem-41}\\
	&R_t=\int_{\Omega}\frac{g(\phi)}{2\sqrt{E_{1}(\phi)}}\phi_{t}\ d\x,\label{problem-51}\\
	&T_{t}=D \Delta T+\frac{R(t)}{\sqrt{E_{1}(\phi)}}  K h^{\prime}(\phi) M(\phi)\mu, \label{problem-61}
\end{align}
\ese
where $M(\phi)=\frac{1}{\varrho(\phi)}$,  and
\begin{equation}
g(\phi)=-\nabla \cdot\left(\left(\kappa^{2}(\nabla \phi)-S_{1}\right) \nabla \phi+\kappa(\nabla \phi)|\nabla \phi|^{2} \boldsymbol{H}(\phi)\right)+\frac{1}{\varepsilon^{2}}\left(f(\phi)-S_{2} \phi\right).\nonumber
\end{equation}
Obviously, the equation \eqref{problem-51} can be obtained by taking the time derivative of the auxiliary variable $R(t)$.
The initial conditions for \eqref{sche11} take
\bex
\phi|_{t=0} =\phi_{0}, \ \ T|_{t=0}=T_{0}, \ \ 
R|_{t=0} =\sqrt{E_{1}(\phi_0)}.
\eex

\begin{rem}
	To separate the computation of different unknown functions, 
	one may think about the splitting method used in \cite{ZCY19novel} 
	or explicit treatment of $\phi_t$ as in \cite{Y21anovel}.
	However the former is unlikely to lead to a second-order scheme, 
	while the latter may result in expensive 
	four-layer computation \cite{Y21anovel}. To construct more efficient second-order scheme, our idea here is to split 
	\eqref{problem-3} into \eqref{problem-31} and \eqref{problem-41}. 
	Another notable idea is to replace $\phi_{t}$ in \eqref{problem-4} by using the equation
	\eqref{problem-31}, resulting in an equivalent equation, i.e., 
	\eqref{problem-61}. 
	Although, and obviously, \eqref{problem-3}-\eqref{problem-4} and \eqref{sche11}
	is strictly equivalent to each other at the continuous level, we will see in what follows that the reformulation \eqref{sche11} facilitates construction of  
	decoupled, stable, higher order convergent, and cheaper time stepping schemes.
\end{rem}

Since \eqref{problem-3}-\eqref{problem-4} and \eqref{sche11}
are equivalent, the latter obviously satisfies  
the same energy dissipative law as \eqref{energylaw}. However, to better understand the discrete 
energy dissipative law, it is desirable to derive an alternative form of the energy law involving the auxiliary variables 
for \eqref{sche11}. To this end,  
we take the $L^{2}$ inner product of $\eqref{problem-31}$ with $-\mu$, 
$\eqref{problem-41}$ with $\phi_{t}$, \eqref{problem-51} with $-2 R$, and 
\eqref{problem-61} with $\frac{\lambda}{\varepsilon K} T$, 
then we perform integration by parts and sum up all equalities to get
\begin{equation}\label{E31}
\frac{d}{d t} E(\phi, R, T)=-\left\|\sqrt{\varrho(\phi)} \phi_{t}\right\|^{2}-\frac{\lambda D}{\varepsilon K}\|\nabla T\|^{2}\le 0,
\end{equation}
where $E(\phi, R, T)$ is defined in \eqref{E3}.

Now we are in a position to construct and analyze our schemes for 
the anisotropic phase-field dendritic crystal growth model \eqref{problem-3}-\eqref{bc}.
To better follow the main argument, let's start with a first-order scheme.

\section {A first order scheme and stability analysis}

Let $\tau>0$ be the time step size, $t^{n}=n \tau, 0 \leq n \leq N$, $T=N \tau$. 
We propose the following scheme:
assuming $\{\phi^{n},T^{n}, R^{n}\}$ are known, $\{\phi^{n+1},T^{n+1}, R^{n+1}\}$ is computed 
by solving:
\bse\label{sche1}
\bal
&\frac{\phi^{n+1}-\phi^{n}}{\tau}=M(\phi^n)\big(\mu^{n+1}-\frac{S_{3}}{\varepsilon^2}(\phi^{n+1}-\phi^n)+S_{4}(\Delta\phi^{n+1}-\Delta\phi^n)\big),\label{problem-311}\\
&\mu^{n+1}=-{\xi}^{n+1}g(\phi^n)+S_{1} \Delta \phi^{n+1}-\frac{S_{2}}{\varepsilon^{2}} \phi^{n+1}-{\xi}^{n+1}\frac{\lambda}{\varepsilon} h^{\prime}(\phi^n) T^{n},\label{problem-411}\\
&\frac{R^{n+1}-R^n}{\tau}= \frac{1}{2\sqrt{E_{1}(\phi^n)}}\left\{ \Big(g(\phi^n),\frac{\phi^{n+1}-\phi^n}{\tau}\Big)-\Big(\frac{\lambda}{\varepsilon}h^{\prime}(\phi^n) M(\phi^n),\mu^n T^{n+1}-\mu^{n+1}T^n\Big)\right . \nonumber\\
&\quad\quad\quad\quad\quad\quad\left. -\Big(\frac{\lambda}{\varepsilon}h^{\prime}(\phi^n)M(\phi^n)T^n ,\frac{S_3}{\varepsilon^2}(\phi^{n+1}-\phi^n)-S_4\Delta(\phi^{n+1}-\phi^n)\Big)\right\},\label{problem-511}\\
&\frac{T^{n+1}-T^n}{\tau}=D\Delta T^{n+1}+{\xi}^{n+1}Kh^{\prime}(\phi^n)M(\phi^n)\mu^n, \label{problem-611}\\
&\left. \frac{\partial\phi^{n+1}}{\partial \n}\right|_{\partial \Omega}=0, \quad \left. 
\frac{\partial T^{n+1}}{\partial \n}\right|_{\partial \Omega}=0, \label{problem-711}
\end{align}
\ese
where ${\xi}^{n+1}=\frac{ R^{n+1}}{\sqrt{E_{1}({\phi}^{n})}}$, 
$S_{3}$ and $S_{4}$ are two extra positive stabilization parameters.

Before carrying out the stability analysis, 
the scheme \eqref{sche1} is worthy of some explanation. 
First, we notice that the coupling terms in both the phase field equation and 
the temperature equation are treated explicitly. This is for ease of calculation.
The implicit treatment of the coupling terms 
in the auxiliary variable equation, i.e., eq.\eqref{problem-511}, may make 
the implementation difficult. However, as we will see in the next section, 
the extra terms added to the scheme, i.e., terms involving the parameters
$S_i$, play a dual role. On one side, some of extra terms are useful in  
decoupling the calculation of different unknowns. 
On the other side, some other extra terms help in 
enhancing the stability. For example, the term 
$\frac{S_{3}}{\varepsilon^{2}}\left(\phi^{n+1}- \phi^{n}\right)$ is 
used to balance the explicit treatment of $\frac{1}{\varepsilon^{2}} f(\phi)$ 
in the phase field equation, and the term $S_{4} \Delta\left(\phi^{n+1}- \phi^{n}\right)$ has purpose to balance the explicit treatment of the gradient term.
The last point we want to emphasize is that the extra terms  
have the same order as the approximation to the time derivatives, thus 
do not affect the overall accuracy. For example, we can check that the term
\bex
&&\frac{1}{2\sqrt{E_{1}(\phi^n)}}\left\{\left(\frac{\lambda}{\varepsilon}h^{\prime}(\phi^n) M(\phi^n),\mu^n T^{n+1}-\mu^{n+1}T^n-\frac{S_3}{\varepsilon^2}T^n(\phi^{n+1}-\phi^n)+S_4T^n\Delta(\phi^{n+1}-\phi^n)\right)\right\}\nonumber
\eex
in \eqref{problem-511} is of order $O(\tau)$. Therefore, formally the convergence of the 
scheme \eqref{sche1} is first order.

\subsection{Stability analysis}

In the following theorem, we establish the stability result for the scheme \eqref{sche1}. That is, we prove
that a discrete ``energy" decays in time, and
consequently the discrete solution remains bounded during the time stepping.

\begin{thm}\label{Stability1}
Let $\{\phi^{n},T^{n}, R^{n}\}$ be the solution of the discrete problem \eqref{sche1}.  Then the following discrete energy law holds:
\begin{equation}\label{energy1}
E^{n+1}-E^n=-Q^{n+1}-\tau\Big(||\sqrt{\varrho(\phi^n)}\frac{\phi^{n+1}-\phi^n}{\tau}||^2+\frac{\lambda D}{\varepsilon K}||\nabla T^{n+1}||^2\Big),
\end{equation}
where $E^{n}$ is defined by
\bex
E^{n}=\frac{S_{1}}{2}||\nabla\phi^n||^2+\frac{S_{2}}{2\varepsilon^2}||\phi^n||^2+\frac{\lambda}{2\varepsilon K} ||T^n||^2
+|R^n|^2,
\eex
$Q^{n+1}$ is given by
\bex
&&Q^{n+1}=\frac{S_1+S_4}{2}||\nabla\phi^{n+1}-\nabla \phi^n||^2+\frac{S_2+S_3}{2\varepsilon^2 }||\phi^{n+1}-\phi^n||^2+ |R^{n+1}-R^n|^2.
\eex
\end{thm}
\begin{proofed}
By taking the inner product of \eqref{problem-311} with $ \frac{2(\phi^{n+1}-\phi^n)}{M(\phi^n)}$, and 
\eqref{problem-411} with $2{(\phi^{n+1}-\phi^n)}$, then summing up the resulting equations, we obtain
\begin{eqnarray}
&&\nonumber\frac{2}{\tau }\left\|\frac{\phi^{n+1}-\phi^n}{\sqrt{M(\phi^n)}}\right\|^2+\frac{2S_3}{\varepsilon^2 }\left\|{\phi^{n+1}-\phi^n}\right\|^2+{2 S_4}\left\|{\nabla(\phi^{n+1}-\phi^n)}\right\|^2+2({\xi}^{n+1}g(\phi^n),\phi^{n+1}-\phi^n) \\
&&+S_1(||\nabla\phi^{n+1}||^2-||\nabla\phi^n||^2+||\nabla(\phi^{n+1}-\phi^n)||^2)+\frac{S_2}{\varepsilon^2 }(||\phi^{n+1}||^2-||\phi^n||^2+||(\phi^{n+1}-\phi^n)||^2)\nonumber\\
&&+2{\xi}^{n+1}\Big(\frac{\lambda}{\varepsilon}h^{\prime}(\phi^n)T^n ,{\phi^{n+1}-\phi^n}\Big)=0.\label{thm141}
\end{eqnarray}
Using \eqref{problem-311}, we rewrite \eqref{problem-511} as follows
\bex
&&\frac{R^{n+1}-R^n}{\tau}= \frac{1}{2\sqrt{E_{1}(\phi^n)}}\left\{ \Big(g(\phi^n),\frac{\phi^{n+1}-\phi^n}{\tau}\Big)-\Big(\frac{\lambda}{\varepsilon}h^{\prime}(\phi^n) M(\phi^n),\mu^n T^{n+1}\Big)\right. \nonumber\\
&&\left.\quad\quad\quad\quad\quad\quad+\Big(\frac{\lambda}{\varepsilon}h^{\prime}(\phi^n)T^n ,\frac{\phi^{n+1}-\phi^n}{\tau}\Big)\right\}.
\eex
Multiplying \eqref{problem-511} with $4\tau R^{n+1}$, we deduce
\begin{eqnarray}
&&2(|R^{n+1}|^2-|R^n|^2+|R^{n+1}-R^n|^2)-2{\xi}^{n+1}\Big\{( g(\phi^{n}),\phi^{n+1}-{\phi}^n)-\tau\Big(\frac{\lambda}{\varepsilon}h^{\prime}(\phi^n) M(\phi^n),\mu^n T^{n+1}\Big) \nonumber\\ && +\Big(\frac{\lambda}{\varepsilon}h^{\prime}(\phi^n)T^n ,{\phi^{n+1}-\phi^n}\Big)\Big\}=0. \label{problem-3611}
\end{eqnarray}
Furthermore, by taking the inner product of  \eqref{problem-611} with $\frac{2\tau \lambda}{\varepsilon K}T^{n+1}$, we obtain
\begin{eqnarray}
&&\frac{\lambda}{\varepsilon K}(||T^{n+1}||^2-||T^n||^2+||T^{n+1}-T^n||^2)+\frac{2\tau\lambda D}{\varepsilon K}||\nabla T^{n+1}||^2\nonumber\\&&-\frac{2\tau\lambda{\xi}^{n+1}}{\varepsilon }\left(h^{\prime}(\phi^n)M(\phi^n)\mu^n, \ T^{n+1}\right)=0.\label{thm1411}
\end{eqnarray}
Finally, the desired result \eqref{energy1} follows from summing up  \eqref{thm141}, \eqref{problem-3611}, and \eqref{thm1411}.
\end{proofed}

\begin{rem}
It is seen from the proof of Theorem \ref{Stability1} that the 
$S_3$ and $S_4$-terms introduced 
in \eqref{problem-311} plays no role in stabilizing the scheme. In fact, it follows from \eqref{energy1}
that the discrete energy remains dissipative even if $S_3=S_4=0$. 
\end{rem}

\subsection{Implementation technique}

It is clear that the efficiency of the scheme depends on whether it can be implemented in an efficient way. 
Besides the provable unconditional stability, 
we will show in this subsection that the proposed scheme can be equivalently reformulated 
into a set of linear elliptic equations, which can be easily solved. 

Noticing that ${\xi}^{n+1}$ is only a scalar variable, we decompose the solution 
$\{\phi^{n+1},\mu^{n+1},T^{n+1}\}$ into the linear combinations as follows: 
\begin{equation}
\label{Com1}
\begin{cases}
\phi^{n+1}=\phi_{1}^{n+1}+{\xi}^{n+1}\phi_{2}^{n+1},
\\
\mu^{n+1}=\mu_{1}^{n+1}+{\xi}^{n+1}\mu_{2}^{n+1},
\\T^{n+1}=T_{1}^{n+1}+{\xi}^{n+1}T_{2}^{n+1}.
\end{cases}
\end{equation}
We impose for the components $\phi^{n+1}_{i}$ and $T_i^{n+1}, i=1,2$, 
the same boundary condition 
as for $\phi^{n+1}$ and $T^{n+1}$, respectively.
Then the equations \eqref{problem-311},  \eqref{problem-411}, and \eqref{problem-611} can be rewritten as
\be\label{mu1e}
\left\{
\ba{l}
\dps\frac{\phi_{1}^{n+1}-\phi^n}{\tau}=M(\phi^n)\big(\mu^{n+1}_1-\frac{S_{3}}{\varepsilon^2}(\phi^{n+1}_1-\phi^n)+S_{4}(\Delta\phi^{n+1}_1-\Delta\phi^n)\big), \\[3mm]
\dps \mu^{n+1}_1=S_{1} \Delta \phi^{n+1}_1-\frac{S_{2}}{\varepsilon^{2}} \phi^{n+1}_1.
\ea
\right.
\ee
\be\label{mu2e}
\left\{
\ba{l}
\dps\frac{\phi_{2}^{n+1}}{\tau}=M(\phi^n)\big(\mu^{n+1}_2-\frac{S_{3}}{\varepsilon^2}\phi^{n+1}_2+S_{4}\Delta\phi^{n+1}_2\big), \\[3mm]
\dps \mu^{n+1}_2=-g(\phi^n)+S_{1} \Delta \phi^{n+1}_2-\frac{S_{2}}{\varepsilon^{2}} \phi_2^{n+1}-\frac{\lambda}{\varepsilon} h^{\prime}(\phi^n) T^{n}.
\ea
\right.
\ee
\begin{numcases}{}
\frac{T^{n+1}_1-T^n}{\tau}=D\Delta T^{n+1}_1,\label{T1e}\qquad\qquad\qquad\qquad\\
\frac{T^{n+1}_2}{\tau}=D\Delta T^{n+1}_2+Kh^{\prime}(\phi^n)M(\phi^n)\mu^n.\label{T2e}\qquad\qquad\qquad
\end{numcases}
Note that $R^{n+1}={\xi}^{n+1}\sqrt{E_{1}(\phi^n)}$, it follows from \eqref{problem-511}:
\begin{eqnarray}
{\xi}^{n+1}A_1^{n+1}=A_2^{n+1},\label{xiv}
\end{eqnarray}
where
\begin{eqnarray}
&&A_1^{n+1}=2E_1(\phi^n)-(g(\phi^n),\phi_2^{n+1})+\frac{\tau\lambda}{\varepsilon}\big(h^{\prime}(\phi^n)M(\phi^n),\mu^nT_2^{n+1}-T^n\mu_2^{n+1}\big)\nonumber\\
&&\qquad\quad\ \ +\frac{\tau\lambda}{\varepsilon  }\big(h^{\prime}(\phi^n)M(\phi^n)T^n,\frac{S_3}{\varepsilon^2}\phi^{n+1}_2-S_4\Delta\phi_2^{n+1}\big),\label{An1}
\end{eqnarray}
\begin{eqnarray}
&&A_2^{n+1}=2\sqrt{E_1(\phi^n)}R^n+(g(\phi^n),\phi_1^{n+1}-\phi^n)-\frac{\tau\lambda}{\varepsilon}\big(h^{\prime}(\phi^n)M(\phi^n),\mu^nT_1^{n+1}-T^n\mu_1^{n+1}\big)\nonumber\\
&&\qquad\qquad-\frac{\tau\lambda}{\varepsilon  }\big(h^{\prime}(\phi^n)M(\phi^n)T^n,\frac{S_3}{\varepsilon^2}(\phi^{n+1}_1-\phi^n)-S_4\Delta(\phi_1^{n+1}-\phi^n)\big).\label{Bn1}
\end{eqnarray}
To simplify the terms above, we use \eqref{mu2e} to obtain
\begin{eqnarray}
&&\quad-(g(\phi^n),\phi_2^{n+1})-\frac{\tau\lambda}{\varepsilon}\big(h^{\prime}(\phi^n)M(\phi^n),T^n\mu_2^{n+1}\big)+\frac{\tau\lambda}{\varepsilon  }\big(h^{\prime}(\phi^n)M(\phi^n)T^n,\frac{S_3}{\varepsilon^2}\phi^{n+1}_2-S_4\Delta\phi_2^{n+1}\big)\nonumber\\
&&=-(g(\phi^n),\phi_2^{n+1})-\big(\frac{\lambda}{\varepsilon  }h^{\prime}(\phi^n)T^n,\phi^{n+1}_2\big)\nonumber\\
&&=\big(\mu^{n+1}_2-S_{1} \Delta \phi^{n+1}_2+\frac{S_{2}}{\varepsilon^{2}} \phi_2^{n+1},\phi^{n+1}_2\big)\nonumber\\
&&=\big( \frac{\phi_2^{n+1}}{\tau M(\phi^n)}+\frac{S_3}{\varepsilon^2}\phi_2^{n+1}-{S_4}\Delta\phi_2^{n+1}-S_1\Delta\phi_2^{n+1}+\frac{S_2}{\varepsilon^2}\phi_2^{n+1},\phi^{n+1}_2 \big)\nonumber\\
&&=\left\|\sqrt{\frac{\varepsilon^2{\varrho(\phi^n)}+\tau(S_2+S_3)}{\tau \varepsilon^2}}\phi_2^{n+1}\right\|^2+\left\|\sqrt{{S_1+S_4}}\nabla\phi_2^{n+1}\right\|^2.\label{A11}
\end{eqnarray}
Then taking the inner product of \eqref{T2e} with $\frac{\tau\lambda}{\varepsilon K}T_2^{n+1}$, we have 
\begin{eqnarray}
&&\quad\frac{\tau\lambda}{\varepsilon}\big(h^{\prime}(\phi^n)M(\phi^n),\mu^nT_2^{n+1}\big)=\frac{\lambda}{\varepsilon K}\big(\left\|T_2^{n+1}\right\|^2+\tau D\left\|\nabla T_2^{n+1}\right\|^2\big).\label{A12}
\end{eqnarray}
Summing up  \eqref{A11} and \eqref{A12},  we get
\begin{eqnarray}
&&A_1^{n+1}=2E_1(\phi^n)+\left\|\sqrt{\frac{\varepsilon^2\varrho(\phi^n)+\tau(S_2+S_3)}{\tau \varepsilon^2}}\phi_2^{n+1}\right\|^2+\left\|\sqrt{{S_1+S_4}}\nabla\phi_2^{n+1}\right\|^2
\nonumber\\
&&\qquad\quad+\frac{\lambda}{\varepsilon K}\big(\left\|T_2^{n+1}\right\|^2
+\tau D\left\|\nabla T_2^{n+1}\right\|^2)>0.\label{An1r}
\end{eqnarray}
Similarly, we can simplify $A_2^{n+1}$ as follows 
\begin{eqnarray}
&&A_2^{n+1}=2\sqrt{E_1(\phi^n)}R^n-\big(\mu_2^{n+1}-S_1\Delta\phi^{n+1}_2+\frac{S_2}{\varepsilon^2}\phi_2^{n+1},\phi_1^{n+1}-\phi^n\big)\nonumber\\
&&\qquad\quad\ \ +\frac{\lambda}{\varepsilon K}(-T_2^{n+1}+\tau D\Delta T_2^{n+1}, T_1^{n+1}).\label{Bn1r}
\end{eqnarray}
Therefore, ${\xi}^{n+1}$ is uniquely determined by dividing the both sides of \eqref{xiv} by 
$A_1^{n+1}$. 

Based on the above discussion, we arrive at the decoupled algorithm for 
solving the equation set \eqref{problem-311}-\eqref{problem-711} as follows:
Given $\{\phi^{n}, T^{n}, R^{n},\mu^n\}$, 
we update  $\{\phi^{n+1}, T^{n+1}, R^{n+1},\mu^{n+1}\}$ through:

i) Solve the equation \eqref{mu1e} for $\phi_{1}^{n+1},~\mu_1^{n+1}$.

\ \ \ Solve the equation \eqref{mu2e} for $\phi_{2}^{n+1},~\mu^{n+1}_{2}$.

ii) Solve the equation \eqref{T1e} for $T_1^{n+1}$.

\ \ \ \ Solve the equation \eqref{T2e} for $T_2^{n+1}$.

iii) Compute ${\xi}^{n+1}$ and $R^{n+1}$ by using 
\eqref{An1r}, \eqref{Bn1r}, and \eqref{xiv}.

iv) Compute $\phi^{n+1},\ \mu^{n+1}$ and $T^{n+1}$ by \eqref{Com1}.

To summarize, the algorithm involves the solution of two variable coefficient linear elliptic equations, two constant coefficient linear elliptic equations, 
and one algebraic equation. Furthermore, it is a three-layer scheme,
compared to the four-layer scheme proposed in \cite{Y21anovel}. 
In the case of constant mobility, our scheme is reduced to solve four elliptic equations with constant coefficients and one algebraic equation.  

\section{A second order scheme}

\subsection {Construction of the scheme}

For ease of notation we will use 
$\bar\varphi^{n+1}$ to denote $2\varphi^{n}-\varphi^{n-1}$.
The second order scheme we propose reads:
\bse\label{sche2}
\bal
&\frac{3\phi^{n+1}-4\phi^{n}+\phi^{n-1}}{2\tau}
=M(\bar{\phi}^{n+1})\big(\mu^{n+1}-\frac{S_{3}}{\varepsilon^2}(\phi^{n+1}-2\phi^n+\phi^{n-1})+S_{4}\Delta(\phi^{n+1}-2\phi^n+\phi^{n-1})\big),\label{problem-342}\\
&\mu^{n+1}=-{\xi}^{n+1}g(\bar{\phi}^{n+1})+S_{1} \Delta \phi^{n+1}-\frac{S_{2}}{\varepsilon^{2}} \phi^{n+1}-{\xi}^{n+1}\frac{\lambda}{\varepsilon} h^{\prime}(\bar{\phi}^{n+1}) \bar{T}^{n+1},\label{problem-352}\\
&\frac{3R^{n+1}-4R^n+R^{n-1}}{2\tau}= \frac{1}{2\sqrt{E_{1}(\bar{\phi}^{n+1})}}\Big\{ \Big(g(\bar{\phi}^{n+1}),\frac{3\phi^{n+1}-4\phi^{n}+\phi^{n-1}}{2\tau}\Big) \nonumber\\
&\quad\quad\quad\quad\quad\quad-\Big(\frac{\lambda}{\varepsilon}h^{\prime}(\bar{\phi}^{n+1}) M(\bar{\phi}^{n+1}),\bar{\mu}^{n+1} T^{n+1}-\mu^{n+1}\bar {T}^{n+1}\Big)\nonumber\\
&\quad\quad\quad\quad\quad\quad-\Big(\frac{\lambda}{\varepsilon}h^{\prime}(\bar{\phi}^{n+1})M(\bar{\phi}^{n+1})\bar{T}^{n+1} ,\frac{S_3}{\varepsilon^2}(\phi^{n+1}-2\phi^n+\phi^{n-1})-S_4\Delta(\phi^{n+1}-2\phi^n+\phi^{n-1})\Big)\Big\}.\label{problem-362}\\
&\frac{3T^{n+1}-4T^n+T^{n-1}}{2\tau}=D\Delta T^{n+1}+{\xi}^{n+1}Kh^{\prime}(\bar{\phi}^{n+1})M(\bar{\phi}^{n+1})\bar{\mu}^{n+1}, \label{problem-372}\\
&\left. \frac{\partial \phi^{n+1}}{\partial \n}\right|_{\partial \Omega}=0, \quad \left. \frac{\partial T^{n+1}}{\partial \n}\right|_{\partial \Omega}=0,\label{problem-382}
\end{align}
\ese
where 
\be\label{Rb2}
{\xi}^{n+1}=\frac{ R^{n+1}}{\sqrt{E_{1}(\bar{\phi}^{n+1})}}.
\ee
Obviously, to start the calculation the scheme \eqref{sche2} must be accompanied by a suitable 
one-step scheme to compute $\left(\phi^{1}, T^{1}, R^{1}\right)$. This can be done, for example, 
by employing the first step of the scheme \eqref{sche1}.

Intuitively this is a second-order scheme since all involved terms are approximated with second-order precision. Although rigorous proof is not available for the time being, the convergence order of the scheme
will be confirmed through a series of numerical tests. 

\subsection {Stability analysis}
Below we prove the stability of the scheme \eqref{sche2}. 
The stability analysis will make use of the following well-known identities:
\begin{eqnarray}
&&2\varphi^{n+1}({3\varphi^{n+1}-4{\varphi}^n+\varphi^{n-1}})=|\varphi^{n+1}|^2-|\varphi^n|^2 
+ |2\varphi^{n+1}-\varphi^n|^2-|2\varphi^n-\varphi^{n-1}|^2\nonumber\\
&&~\quad\quad\quad\quad\quad\quad\quad\quad\quad\quad\quad\quad\quad+|\varphi^{n+1}-2\varphi^n+\varphi^{n-1}|^2,\label{problem-42}\\ \nonumber
&&(\varphi^{n+1}-2\varphi^{n}+\varphi^{n-1})({3\varphi^{n+1}-4{\varphi}^n+\varphi^{n-1}})
=|\varphi^{n+1}-\varphi^n|^2-|\varphi^n-\varphi^{n-1}|^2\\
&&~\quad\quad\quad\quad\quad\quad\quad\quad\quad\quad\quad\quad\quad\quad\quad\quad\quad\quad\quad\quad+2|\varphi^{n+1}-2\varphi^{n}+\varphi^{n-1}|^2.\label{problem-43}
\end{eqnarray}

\begin{thm}\label{Stability2}
Let $\{\phi^{n},  T^{n}, R^{n}\}$ be the solution of the discrete problem \eqref{sche2}.  
Then for $n\geq 1$ it satisfies
the discrete energy law:
\begin{eqnarray}
&&E^{n+1}-E^n=-Q^{n+1}-\tau\Big(\left\|\sqrt{\varrho(\bar{\phi}^{n+1})}\Big(\frac{3\phi^{n+1}-4\phi^n+\phi^{n-1}}{2\tau}\Big)\right\|^2+\frac{\lambda D}{\varepsilon K}||\nabla T^{n+1}||^2\Big),\label{energy2}
\end{eqnarray}
where $E^{n}$ is defined by
\begin{eqnarray}
&&\!\!\!\!\!\!\!\!E ^{n}=\frac{1}{4}\big[S_1(||\nabla\phi^n||^2+||2\nabla\phi^n-\nabla\phi^{n-1}||^2)+\frac{S_2}{\varepsilon^2}(||\phi^{n}||^2+||2\phi^{n}-\phi^{n-1}||^2)+\frac{2S_3}{\varepsilon^2}||\phi^{n}-\phi^{n-1}||^2\nonumber\\
&&\!\!\!\quad\ +{2S_4}||\nabla(\phi^{n}-\phi^{n-1})||^2+\frac{\lambda}{\varepsilon K}(||T^n||^2+||2T^n-T^{n-1}||^2)+2(|R^n|^2+|2R^n-R^{n-1}|^2)\big],\label{equa2}
\end{eqnarray}
and $Q^{n+1}$ is defined by
\begin{eqnarray}
&&Q^{n+1}= \frac{S_1+2S_4}{4}||\nabla\phi^{n+1}-2\nabla\phi^{n}+\nabla\phi^{n-1}||^2+\frac{S_2+2S_3}{4\varepsilon^2}||\phi^{n+1}-2\phi^n+\phi^{n-1}||^2\nonumber\\
&&\qquad\quad\  +\frac{\lambda}{4\varepsilon K}||T^{n+1}-2T^n+T^{n-1}||^2+\frac{1}{2}|R^{n+1}-2R^n+R^{n-1}|^2.\nonumber
\end{eqnarray}
\end{thm}
\begin{proofed}
First, we take the inner product of  \eqref{problem-342} with $\frac{2}{M(\bar{\phi}^{n+1})}(3\phi^{n+1}-4\phi^n+\phi^{n-1})$, and \eqref{problem-352} with $2{(3\phi^{n+1}-4{\phi}^n+\phi^{n-1})}$.
Then we sum up the resulting equations and use the identities \eqref{problem-42} and \eqref{problem-43} 
to obtain
\begin{eqnarray}
&&\nonumber{4\tau}\left\|\frac{3\phi^{n+1}-4\phi^n+\phi^{n-1}}{2\tau \sqrt{M(\bar{\phi}^{n+1})}}\right\|^2+\frac{2S_3}{\varepsilon^2}(\left\|{\phi^{n+1}-\phi^{n}}\right\|^2-\left\|\phi^{n}-\phi^{n-1}\right\|^2+2\left\|\phi^{n+1}-2\phi^{n}+\phi^{n-1}\right\|^2)\\
&&+{2S_4}\left(\left\|\nabla(\phi^{n+1}-\phi^{n})\right\|^2-\left\|{\nabla(\phi^{n}-\phi^{n-1})}\right\|^2+2\left\|{\nabla(\phi^{n+1}-2\phi^{n}+\phi^{n-1})}\right\|^2\right)\nonumber\\
&&+S_1(||\nabla\phi^{n+1}||^2+||\nabla(2\phi^{n+1}-\phi^n)||^2-||\nabla\phi^{n}||^2-||\nabla(2\phi^{n}-\phi^{n-1})||^2+||\nabla(\phi^{n+1}-2\phi^n+\phi^{n-1})||^2)
\nonumber\\
&&+\frac{S_2}{\varepsilon^2}(||\phi^{n+1}||^2+||2\phi^{n+1}-\phi^n||^2-||\phi^{n}||^2-||2\phi^{n}-\phi^{n-1}||^2+||\phi^{n+1}-2\phi^n+\phi^{n-1}||^2)\nonumber\\
&&+2{\xi}^{n+1}\big(g(\bar{\phi}^{n+1}),{3\phi^{n+1}-4{\phi}^n+\phi^{n-1}}\big)+\frac{4\tau\lambda{\xi}^{n+1}}{\varepsilon}\Big(h^{\prime}(\bar{\phi}^{n+1})\bar{T}^{n+1},\frac{3\phi^{n+1}-4\phi^n+\phi^{n-1}}{2\tau}\Big)=0.\label{thm181} 
\end{eqnarray}
By virtue of \eqref{problem-362} and \eqref{problem-342}, we can rewrite \eqref{problem-362} as 
\bex
&&\frac{3R^{n+1}-4R^n+R^{n-1}}{2\tau}- \frac{1}{2\sqrt{E_{1}(\bar{\phi}^{n+1})}}\Big\{ \Big(g(\bar{\phi}^{n+1}),\frac{3\phi^{n+1}-4\phi^{n}+\phi^{n-1}}{2\tau}\Big) \nonumber\\
&&-\Big(\frac{\lambda}{\varepsilon}h^{\prime}(\bar{\phi}^{n+1}) M(\bar{\phi}^{n+1}),\bar{\mu}^{n+1} T^{n+1}\Big)+\Big(\frac{\lambda}{\varepsilon}h^{\prime}(\bar{\phi}^{n+1})\bar{T}^{n+1} ,\frac{3\phi^{n+1}-4\phi^n+\phi^{n-1}}{2\tau}\Big)\Big\}=0.
\eex
Multiplying the both sides by $8\tau R^{n+1}$, we have
\begin{eqnarray}
&&2(|R^{n+1}|^2+|2R^{n+1}-R^{n}|^2-|R^n|^2-|2R^{n}-R^{n-1}|^2+|R^{n+1}-2R^n+R^{n-1}|^2)\nonumber\\
&&-4\tau{\xi}^{n+1}\Big\{\Big({ g(\bar{\phi}^{n+1})},\frac{{3\phi^{n+1}-4{\phi}^n+\phi^{n-1}}}{2\tau}\Big)
-\Big(\frac{\lambda}{\varepsilon}h^{\prime}(\bar{\phi}^{n+1}) M(\bar{\phi}^{n+1}),\bar{\mu}^{n+1} T^{n+1}\Big)
\nonumber \\
&& +\Big(\frac{\lambda}{\varepsilon}h^{\prime}(\bar{\phi}^{n+1})\bar{T}^{n+1} ,\frac{3\phi^{n+1}-4\phi^n+\phi^{n-1}}{2\tau}\Big)\Big\}=0.\label{thm1911}
\end{eqnarray}
Furthermore, it follows from taking the inner product  of \eqref{problem-372} with $\frac{4\lambda\tau}{\varepsilon K}T^{n+1}$, and using \eqref{problem-42}:
\begin{eqnarray}
&&\frac{\lambda}{\varepsilon K}\left\{(||T^{n+1}||^2+||2T^{n+1}-T^n||^2-||T^{n}||^2-||2T^{n}-T^{n-1}||^2+||T^{n+1}-2T^n+T^{n-1}||^2)\right.
\nonumber\\
&&\left. +{4\tau D}||\nabla T^{n+1}||^2-4\tau{\xi}^{n+1}(Kh^{\prime}(\bar{\phi}^{n+1})M(\bar{\phi}^{n+1})\bar{\mu}^{n+1},T^{n+1}) \right\}=0.\label{thm1811}
\end{eqnarray}
Finally we sum up \eqref{thm181}, \eqref{thm1911}, and  \eqref{thm1811} to conclude.
This completes the proof.
\end{proofed}

The second order scheme \eqref{sche2} can also be implemented in an efficient way. We are not going to describe the implementation details of this scheme since it is very similar to the first order scheme \eqref{sche1} explained in the previous section. 
Nevertheless we give here the algebraic equation to compute the auxiliary 
variable ${\xi}^{n+1}$:
\bex
{\xi}^{n+1}=\frac{A_2^{n+1}}{A_1^{n+1}},
\eex 
where 
\bex
&&A_1^{n+1}=3E_1(\bar{\phi}^{n+1})+\frac{3}{2}\left\|\sqrt{\frac{3\varepsilon^2\varrho(\bar{\phi}^{n+1})+2\tau(S_2+S_3)}{2\tau \varepsilon^2}}\phi_2^{n+1}\right\|^2+\frac{3}{2}\left\|\sqrt{{S_4+S_1}}\nabla\phi_2^{n+1}\right\|^2\\
&&\qquad\quad\ \ +\frac{\lambda}{\varepsilon K}\big(\frac{3}{2}\left\|T_2^{n+1}\right\|^2+\tau D\left\|\nabla T_2^{n+1}\right\|^2\big)>0,
\eex
\bex
&&A_2^{n+1}=\sqrt{E_1(\bar{\phi}^{n+1})}(4R^n-R^{n-1})-\frac{1}{2}\big(\mu_2^{n+1}-S_1\Delta\phi^{n+1}_2+\frac{S_2}{\varepsilon^2}\phi_2^{n+1},3\phi_1^{n+1}-4\phi^n+\phi^{n-1}\big)\\
&&\qquad\quad\ \ + \frac{\lambda}{\varepsilon K}\big(-\frac{3}{2}T_2^{n+1}+\tau D\Delta T_2^{n+1}, T_1^{n+1}\big).
\eex

\section{Numerical experiments}
\setcounter{equation}{0}

In order to illustrate the performance of the proposed numerical method and confirm our analysis results, 
several numerical tests are carried out and presented in this section. 
For the sake of convenience, we consider the numerical examples with 
the square domain $\Omega=(-1,1)^2$ and fourfold anisotropy, i.e., $m=4$
in \eqref{Ka}. We consider constant mobility $\varrho(\phi)$ in the calculation. 
The spatial discretization is a Legendre-Galerkin spectral method. The approximation space for the phase function $\phi$ and the temperature $T$ is $I\!\!P_{N}(\Omega)$, where
$I\!\!P_N(\Omega)$ denotes the space of polynomials of degree $\le N$ at each space direction. 
Since the focus of our numerical tests is the verification of the stability and
convergence order of the time-stepping schemes, we will take 
$N=128$, which is large enough so that the spatial discretization errors are negligible compared with  the temporal one.
\begin{example}\label{exp1}(Accuracy test)
As the first example, we test the convergence order of the schemes \eqref{sche1} and \eqref{sche2}. 
In the case one, we fabricate two forcing functions in \eqref{problem-3} and \eqref{problem-4} so that the exact solution to \eqref{problem-3}-\eqref{bc} is
\be\label{case1}
\mbox{Case-I}\quad\quad\quad \left\{
\begin{array}{rcl}
\phi(x,y,t)=\sin(t)\cos(\pi x)\cos(\pi y),\quad\\
T(x,y,t)=\sin(t)\cos^2( \pi x)\cos^2(\pi y) .
\end{array} \right.
\ee
We set the following parameters:
\begin{align}\label{p-case1}
\left\{
\ba{l}
\varrho=4e3,\  \varepsilon=0.1,\  \sigma=0.05,\  \lambda=0.1,\  D=2.25e-2,\  K=0.01,\\
S_1=0.9,\  S_2=10,\  S_3=S_4=0,\  B=1e4.
\ea
\right.
\end{align}
In the case two, we choose the initial conditions:
\be\label{case2}
\mbox{Case-II}\quad\quad\quad \left\{
\begin{array}{l}
\dps \phi(x,y,0)=\tanh\Big(\frac{r_0-((x-x_0)^2+(y-y_0)^2)}{\varepsilon_0}\Big),\ \ \\[3mm]
T(x,y,0)=-0.5\phi(x,y,0),\qquad\qquad\qquad
\end{array} 
\right. 
\ee
where $r_0=0.25,\ x_0=y_0=0,\ \varepsilon_0=0.1$. 
The model parameters are set as follows:
\begin{align}\label{p-case2}
\left\{
\ba{l}
\varrho=1e3,\  \varepsilon=0.1,\  \sigma=0.05,\  \lambda=1,\    D=5e-2,\  K=0.1,\\
S_1=0.9,\  S_2=10,\  S_3=S_4=0,\  B=5e3.
\ea
\right.
\end{align}
In the latter case the exact solution is unavailable, we will use the numerical solution computed with $\tau=3e-5$ to serve as the exact solution.
In Figure \ref{figcon} we plot the $L^2$ errors in log-log scale of the computed phase and temperature solutions at $t=1$ as functions of the time step size 
$\tau$. As expected, in both cases the obtained numerical convergence rates 
are in a perfect agreement with the claimed orders; i.e., first order for the scheme \eqref{sche1} and second order for the scheme \eqref{sche2}. 
It is worth to mention that the calculation with some 
positive parameters $S_3$ and $S_4$ has produced similar results (not shown here).

\begin{figure}[!ht]
\centering
\subfigure[First order scheme \eqref{sche1} for the Case-I.]{
	\includegraphics[width=7.3cm]{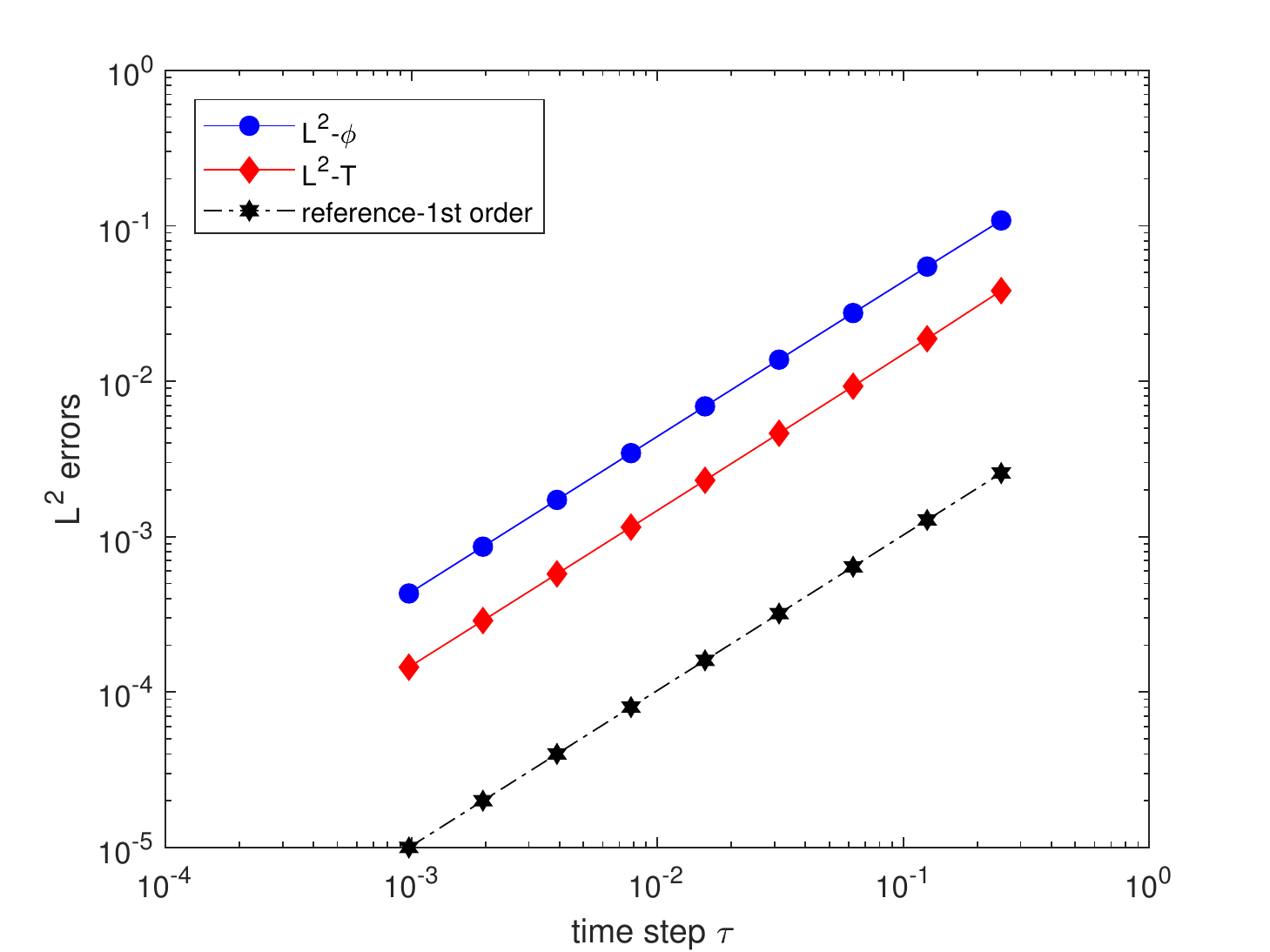}
}
\hspace{-10mm}
\subfigure[First order scheme \eqref{sche1} for the Case-II.]{
	\includegraphics[width=7.3cm]{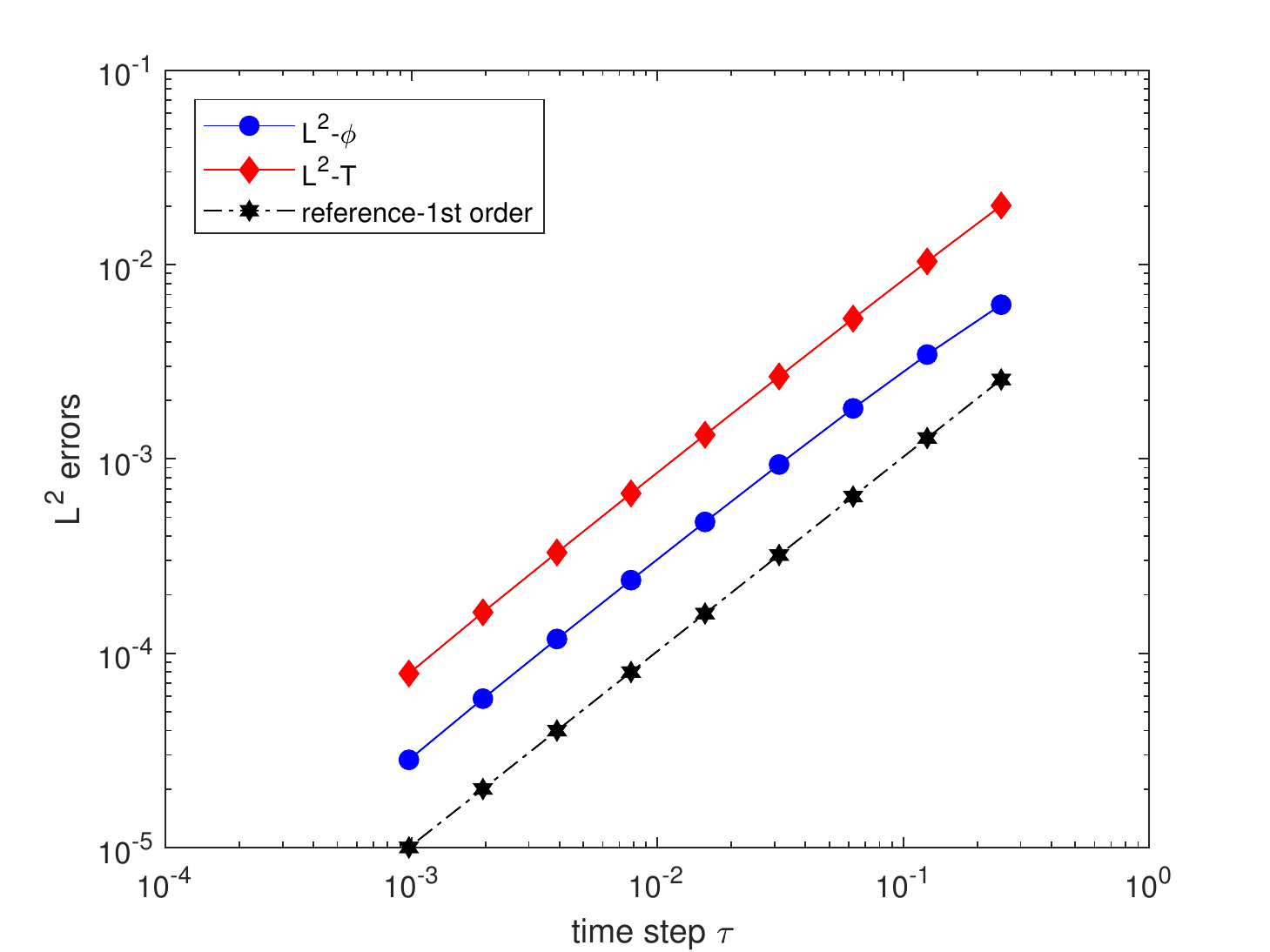}}
\centering
\subfigure[Second order scheme \eqref{sche2} for the Case-I.]{
	\includegraphics[width=7.3cm]{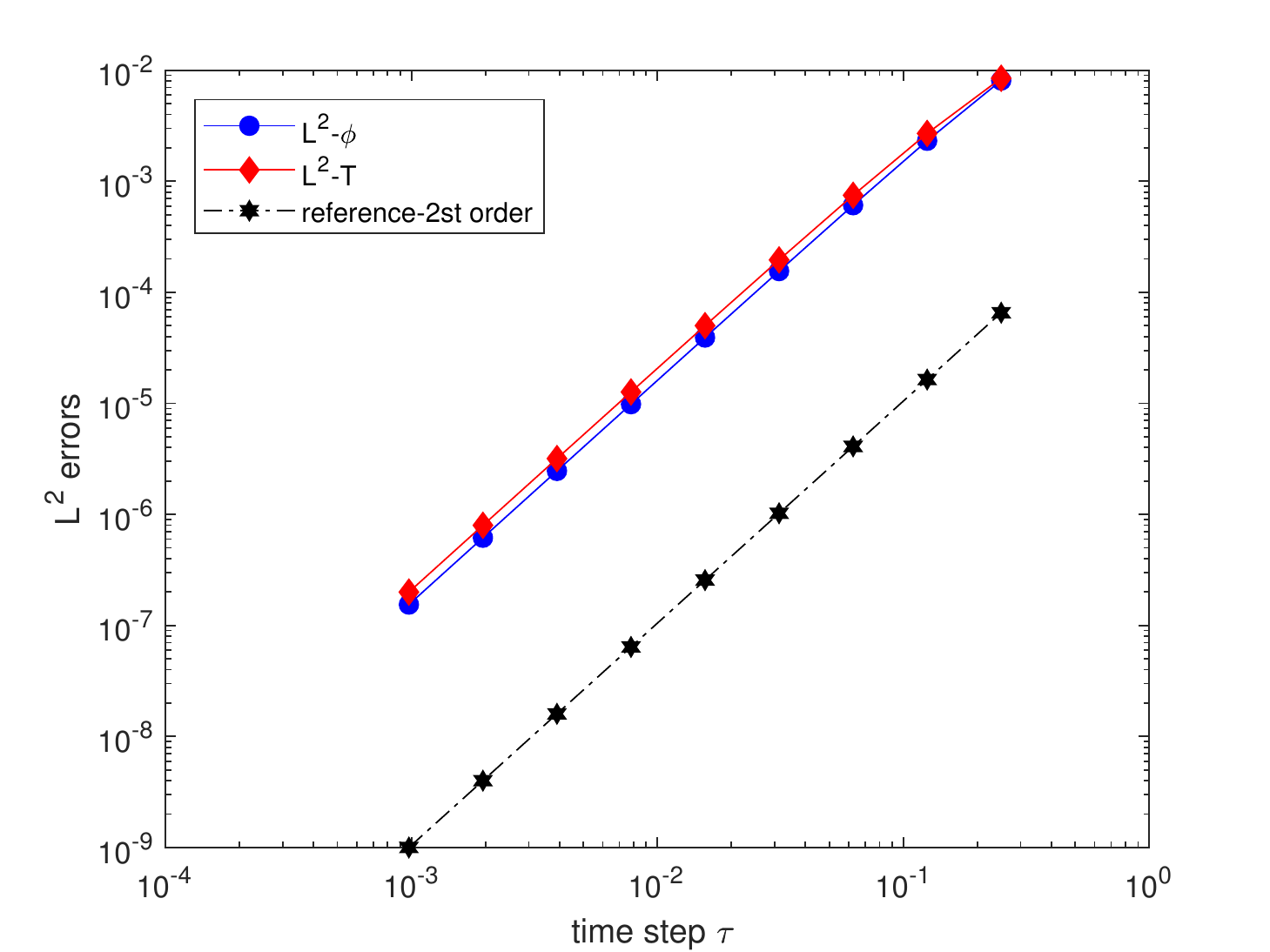}
}
\hspace{-10mm}
\subfigure[Second order scheme \eqref{sche2} for the Case-II.]{
	\includegraphics[width=7.3cm]{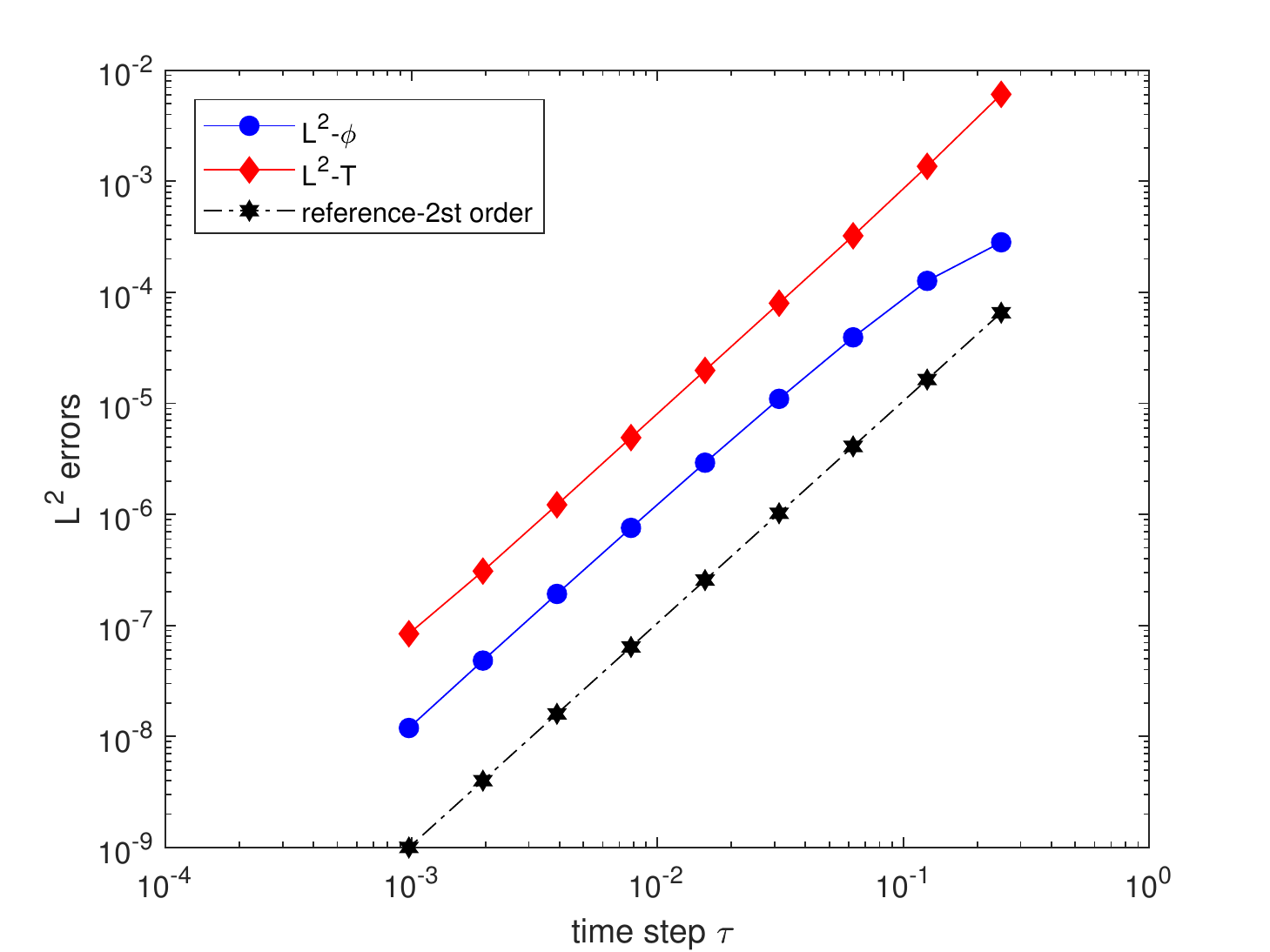}}
\caption{(Example \ref{exp1}) Convergence order of the time-stepping schemes: $L^2$ errors of the phase field function and the temperature as functions of the time step size $\tau$.}
\label{figcon}
\end{figure}
\end{example}

\begin{example}\label{exp2}(Stability test)
To investigate the stability property of the proposed schemes, 
we consider the problem Case-II, which has the initial conditions given in
\eqref{case2}. The parameters used in this test are the same as in 
\eqref{p-case2}. 

We know from \eqref{E3} that $E^n$ defined in \eqref{equa2} 
can be regarded as a discrete version of the original energy functional 
$E(\phi^n, T^n)$ defined in \eqref{E1}. 
According to \eqref{energylaw} and Theorem \ref{Stability2}, both $E(\phi^n, T^n)$ and $E^n$ should be monotonically decreasing with the time step $n\ge 1$. 
Notice that $E^n$ is usually called as modified energy functional, 
which is not necessarily an approximation to the original energy functional 
$E(\phi^n, T^n)$.

\begin{figure}[!ht]
\centering
\subfigure[Evolution of $E^n$]{
	\includegraphics[width=2.08in]{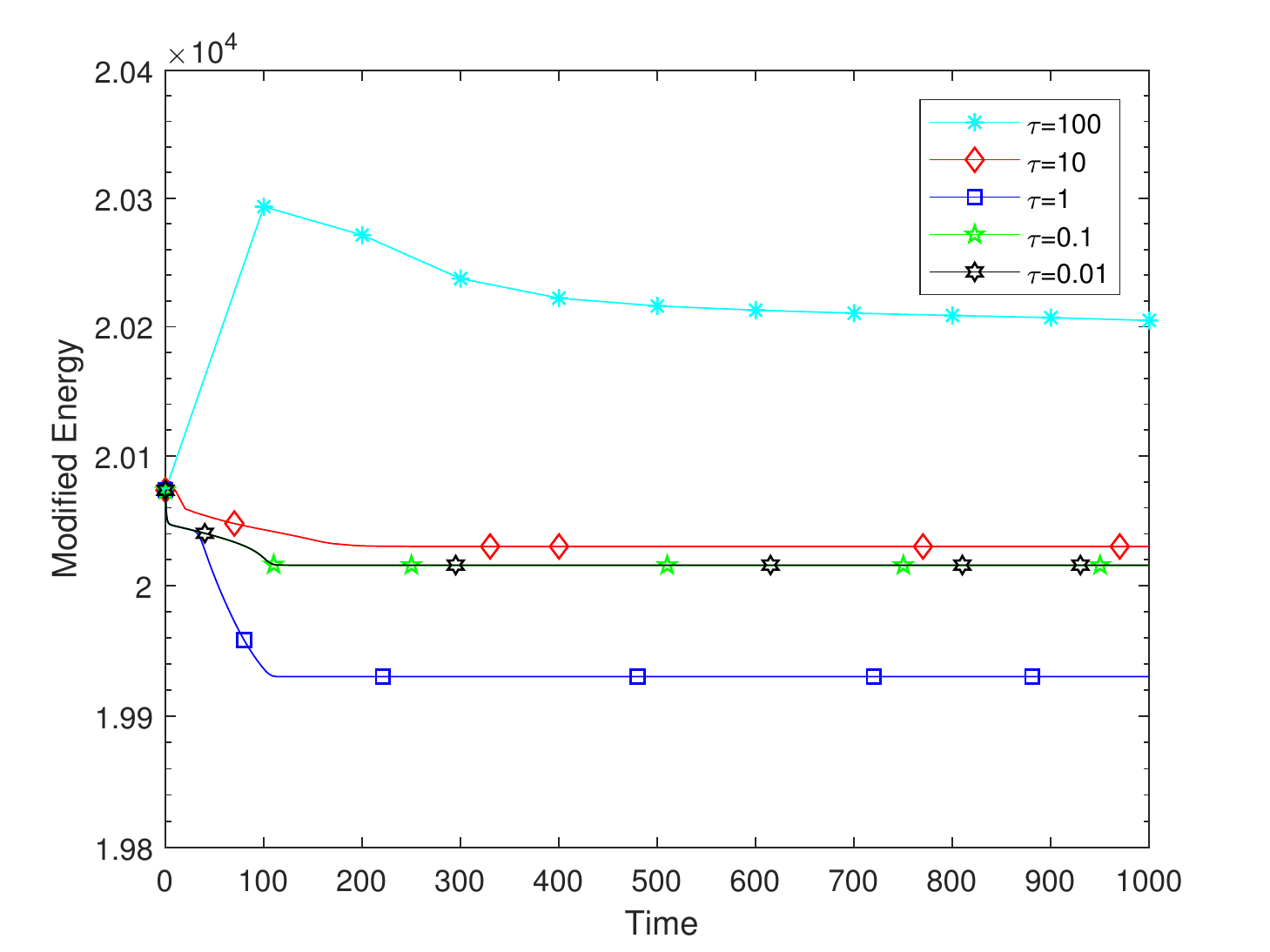}
}
\hspace{-8mm}
\subfigure[Evolution of $E(\phi^n, T^n)$]{
	\includegraphics[width=2.08in]{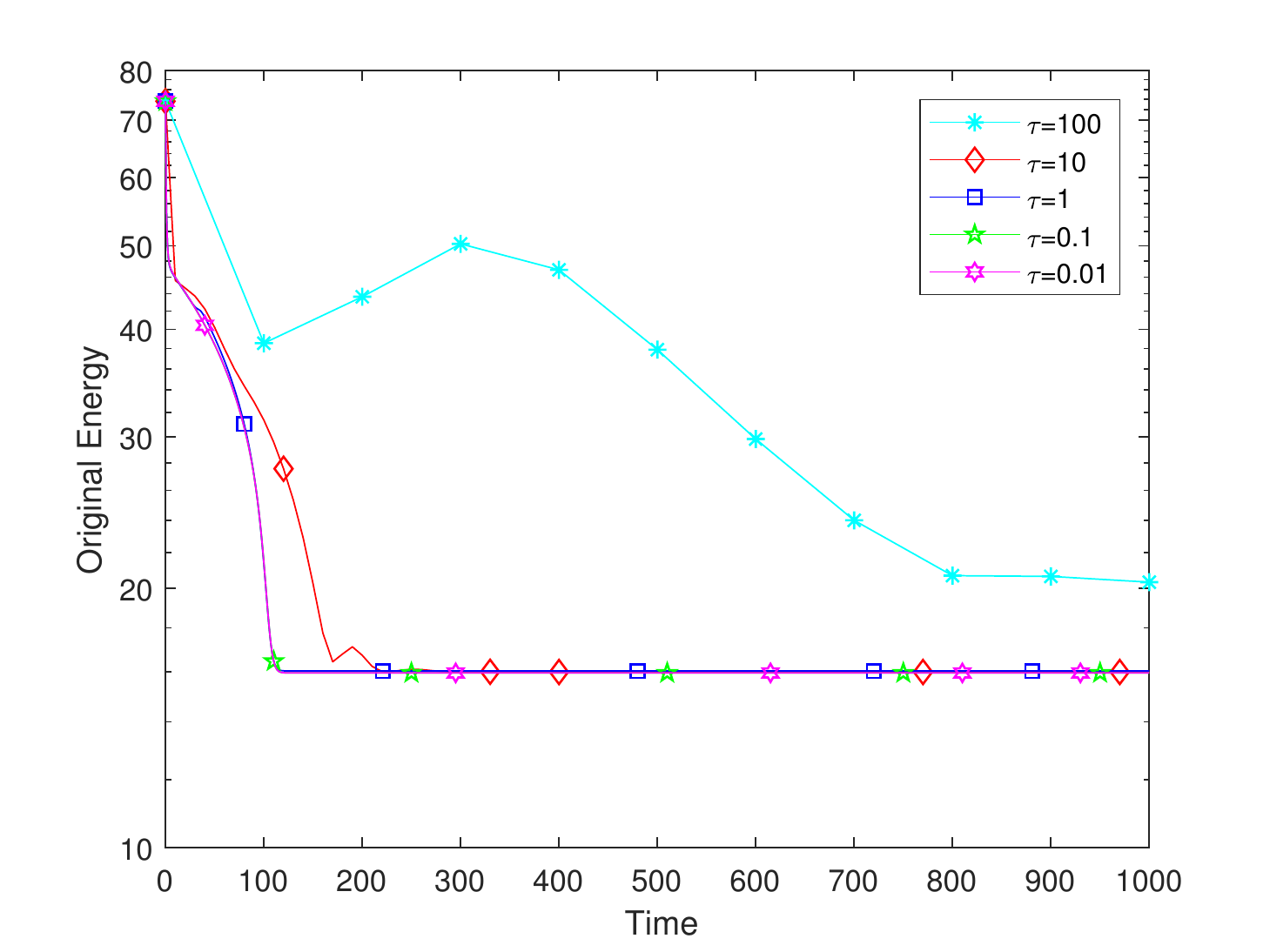}
}
\hspace{-8mm}
\subfigure[Evolution of  ${\xi}^n$]{
	\includegraphics[width=2.08in]{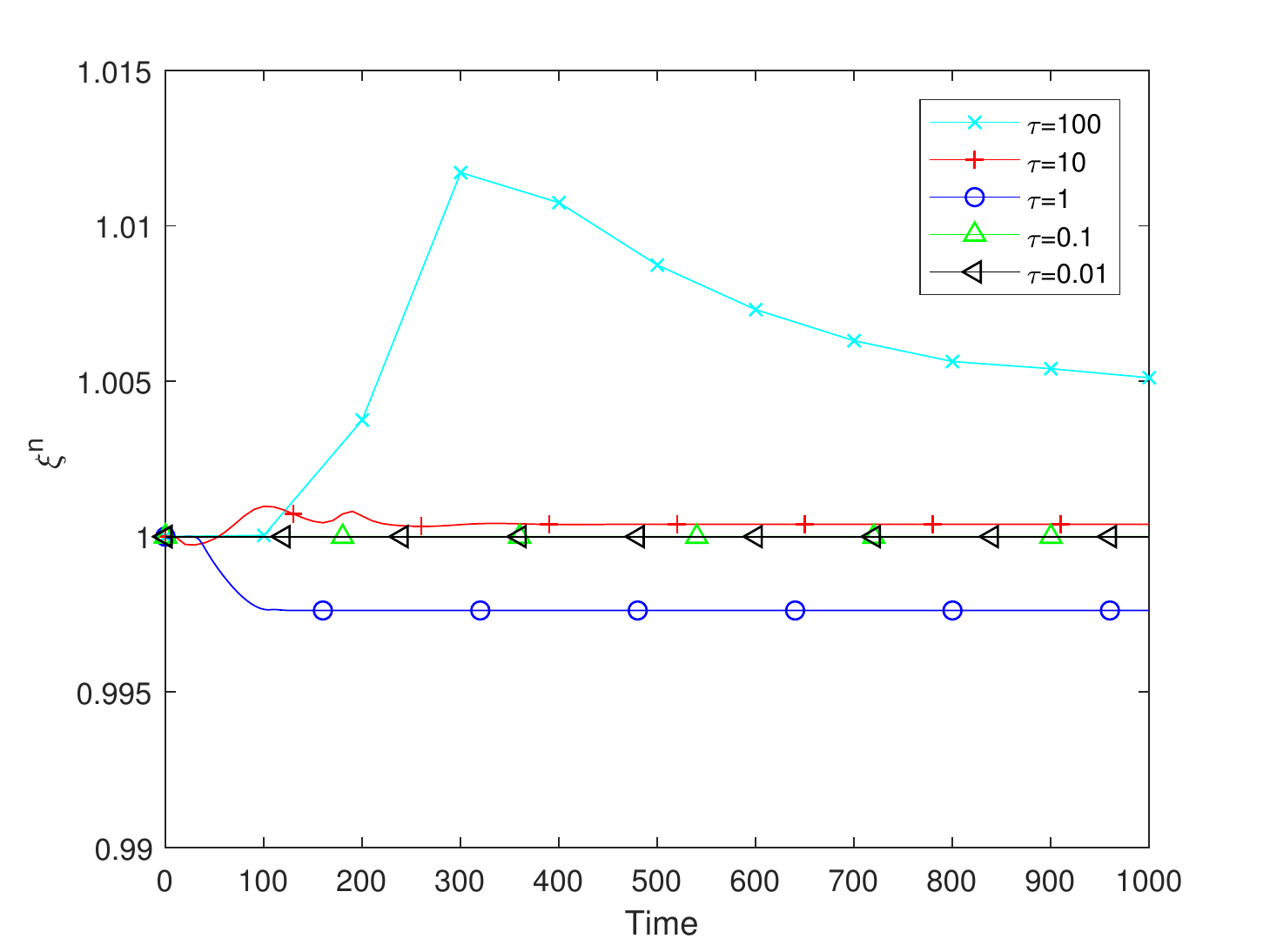}
}
\caption{(Example \ref{exp2})  Time evolution of the modified energy functional $E^n$, the original energy
	$E(\phi^n, T^n)$, and ${\xi}^n$ computed 
	by the scheme \eqref{sche2} using several time step sizes.}
\label{figen1}
\end{figure}
The  modified energy functional $E^n$, original energy functional 
$E(\phi^n, T^n)$, and 
${\xi}^n$ computed by the scheme \eqref{sche2} with different 
time step sizes are presented in Figure \ref{figen1} as functions of 
time.
It is observed in Figure \ref{figen1}(a) that 
the modified energy functional $E^n$
is indeed strictly dissipative in time as predicted by Theorem \ref{Stability2}, even for very large time step sizes $\tau=10, 100$.
As shown in Figure \ref{figen1}(b), 
the original energy $E(\phi^n, T^n)$ is also dissipative for all relatively small time step sizes.
However, the dissipation of $E(\phi^n, T^n)$ loses monotonicity 
during some time period($\tau\geq 10$) for the time step sizes bigger 
than 10. This is most likely caused by imprecise calculation with 
large time step sizes. 
This guess is supported by the observation 
from Figure \ref{figen1}(c), 
in which the time evolution of ${\xi}^n$ is given. Remember that ${\xi}^n$, defined in \eqref{Rb2}, should be close to 1 if the approximation is good enough. The results presented in Figure \ref{figen1}(c) demonstrate good accuracy on ${\xi}^n$ 
when $\tau$ is not very large, say less than $10$. However when 
$\tau=10$, ${\xi}^n$ becomes oscillatory and error becomes visible.
When $\tau$ increases to 100, the computed auxiliary variable 
${\xi}^n$ is not any more close to 1.
This implies that the numerical solution is not accurate enough for 
$\tau\ge 10$, leading to a violation of the monotonic dissipation of the original energy. 
Thus in practice, it is not recommended to use time step sizes too large, although the calculation can always be stable. 

Next test concerns impact of 
the parameters $S_1,\ S_2,\ S_3$, and $S_4$ on the quality of numerical solutions. Since the auxiliary variable is indicative of the quality of the computed solutions, we only report the computed values of ${\xi}^n$. 
Figure \ref{figxi} shows the computed ${\xi}^n$ versus the
time for four parameter sets $\{S_i\}$. It is observed from
Figures \ref{figxi}(a)-(b) that for
$S_1=S_2=S_3=S_4=0$ and $S_1=0.1,\ S_2=4,\ S_3=S_4=0$, 
the computed ${\xi}^n$ is quite inaccurate, i.e., far from the exact 
value 1, even with small time step sizes (0.001 and 0.01 respectively). 
However taking positive $S_3$ and $S_4$, 
say $S_3=S_4=5$, allows recovering the accuracy, as
shown  
in Figures \ref{figxi}(c)-(d).
It is notable that the presence of the $S_3$- and $S_4$-terms allows
stable and accurate calculation even with $\tau=10$. This test clearly 
indicates the benefit of the stabilization terms.
It is also worth to mention that all the cases above 
produced monotonically decreasing energy $E^n$,  
which is consistent with what we have proved in Theorem \ref{Stability2}.

\begin{figure}[!ht]
\centering
\subfigure[$S_1=S_2=S_3=S_4=0$]{
	\includegraphics[width=7.3cm]{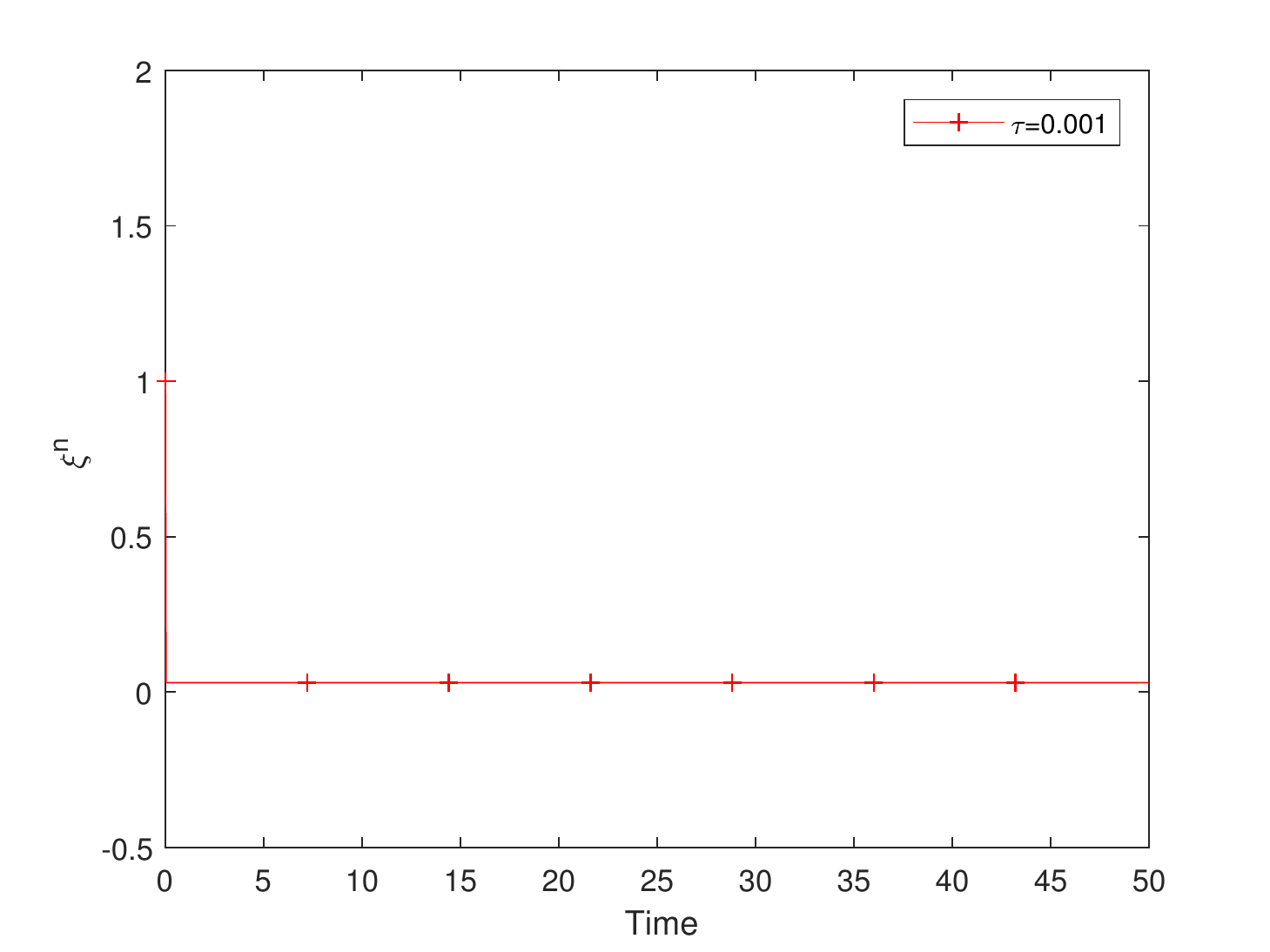}
}
\hspace{-10mm}
\subfigure[$S_1=0.1,\ S_2=4,\ S_3=S_4=0$]{
	\includegraphics[width=7.3cm]{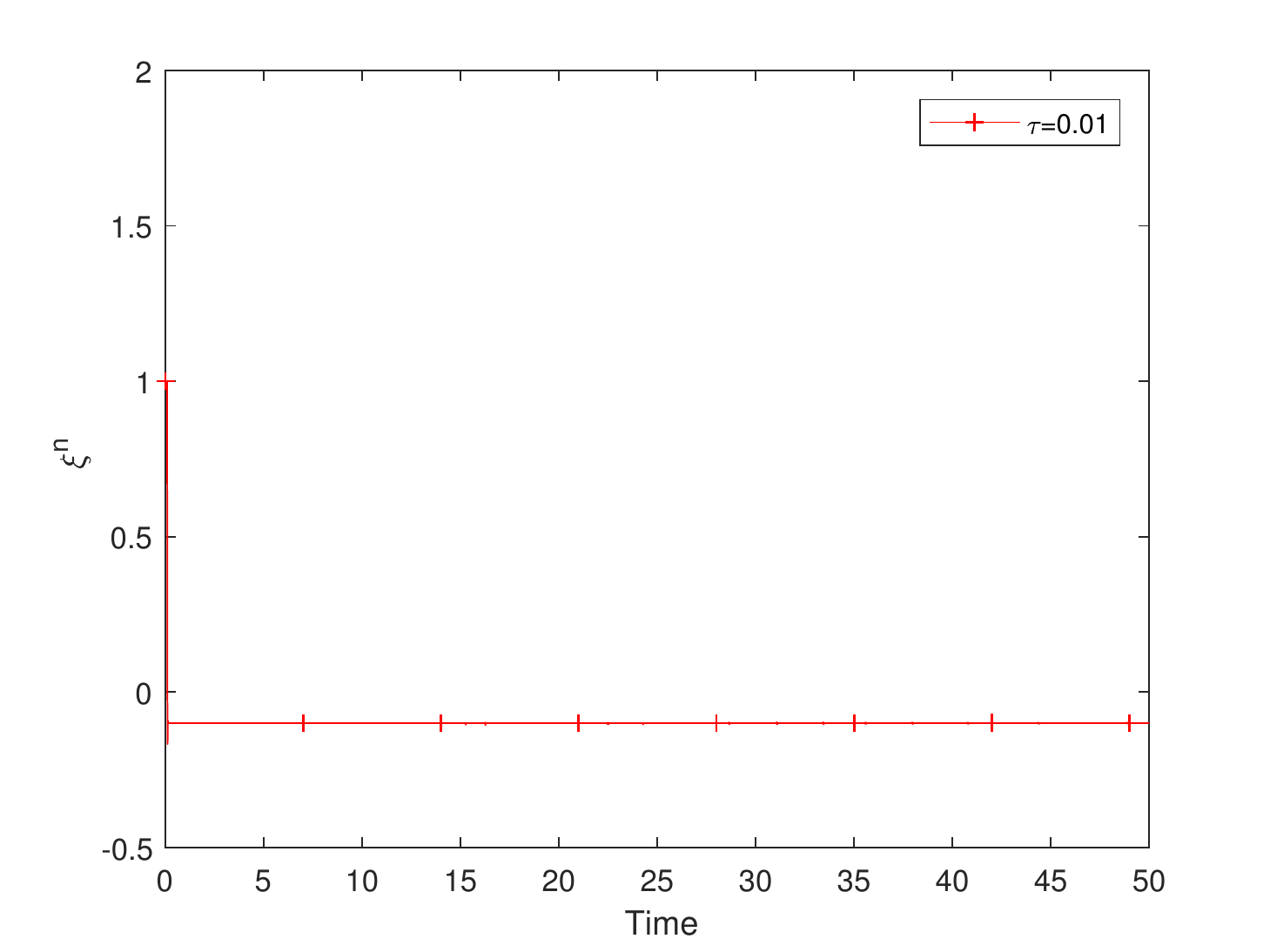}
}
\centering
\subfigure[$S_1=S_2=0,\ S_3=S_4=5$]{
	\includegraphics[width=7.3cm]{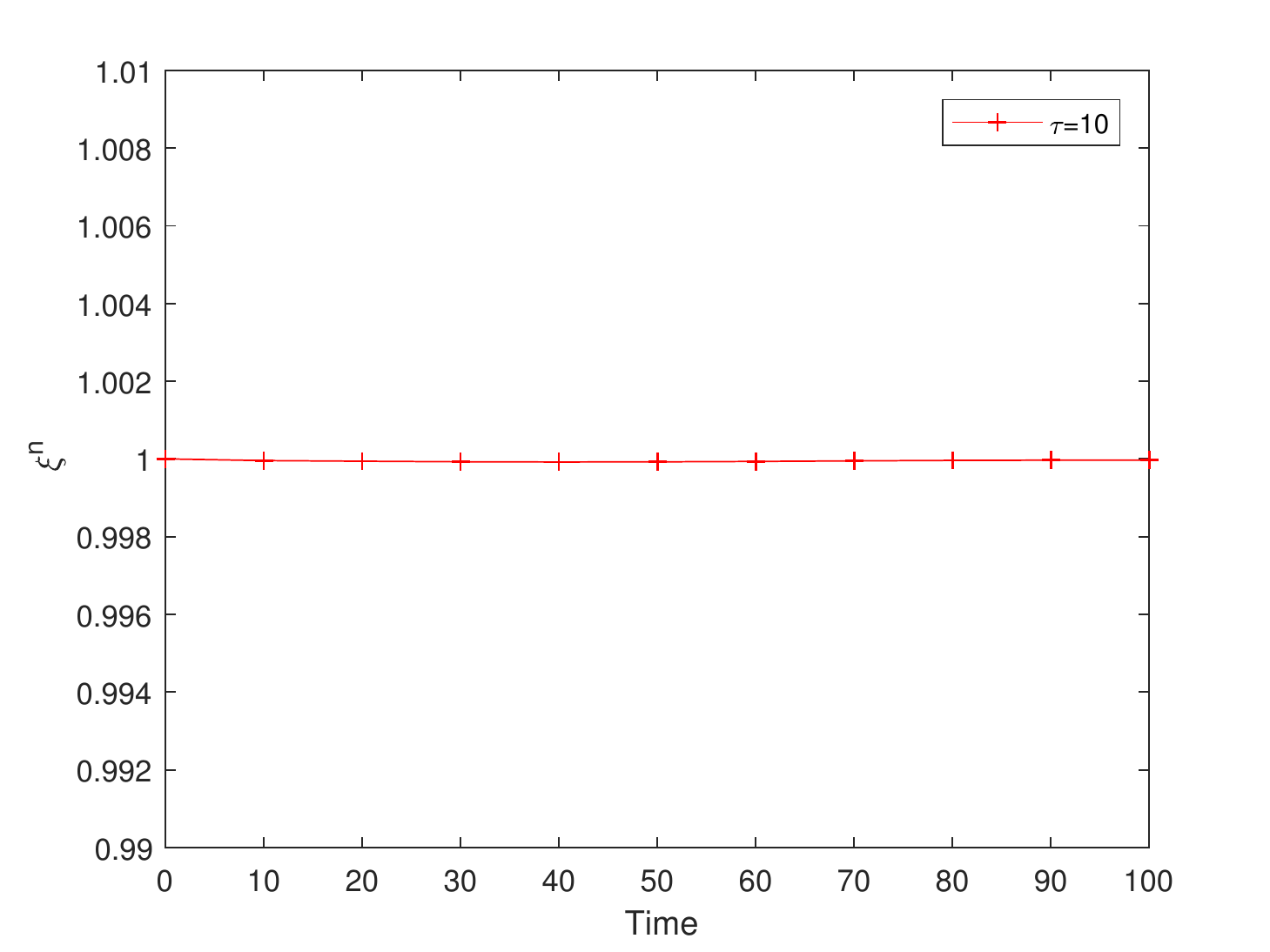}
}
\hspace{-10mm}
\subfigure[$S_1=0.1,\ S_2=4,\ S_3=S_4=5$]{
	\includegraphics[width=7.3cm]{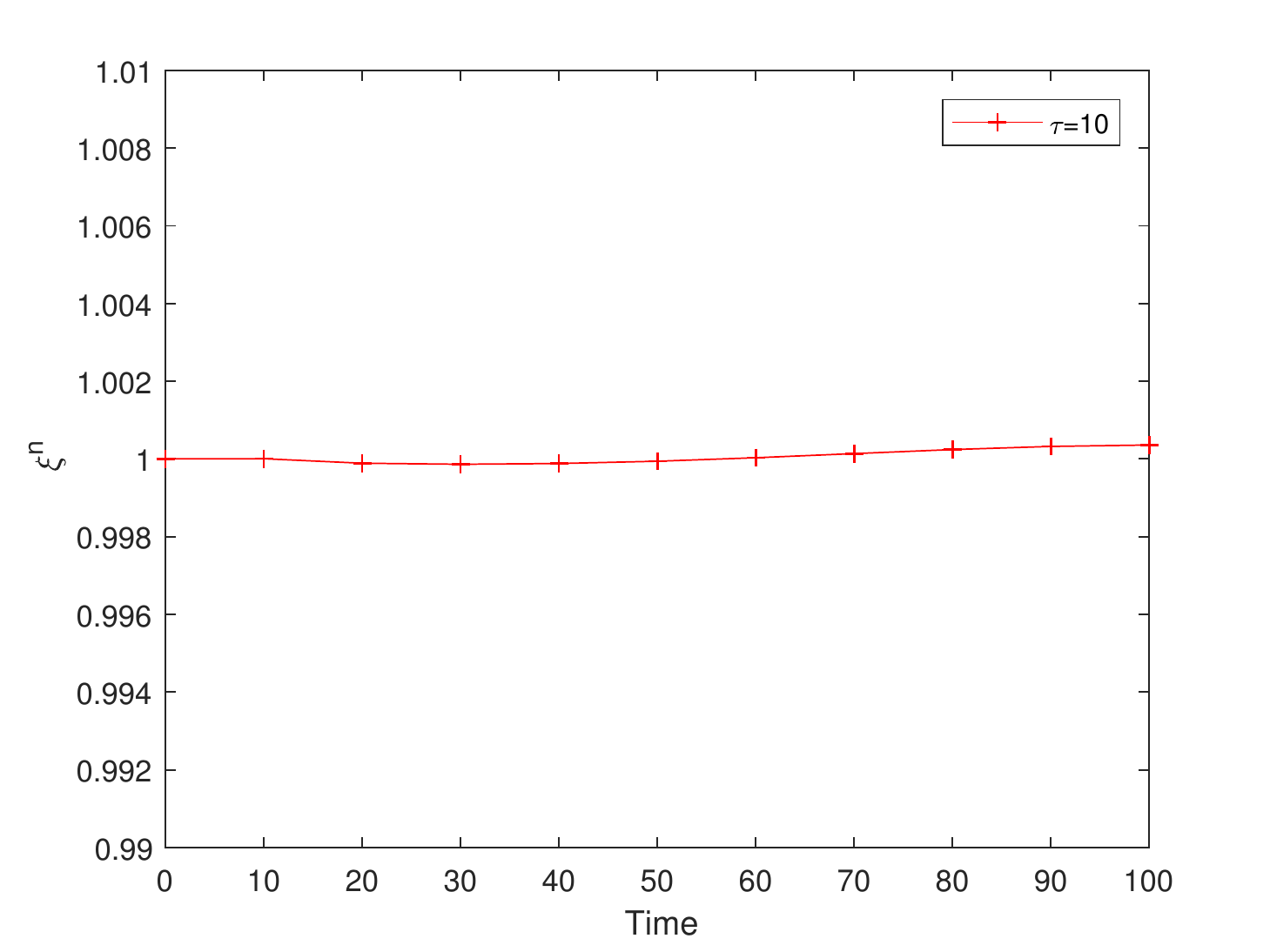}}
\caption{(Example \ref{exp2}) History of ${\xi}^n$ with different chioces of the parameters $S_1,\ S_2,\ S_3,\ S_4$.}
\label{figxi}
\end{figure}
\end{example}

\begin{example}\label{exp3}
(Fourfold anisotropy crystal growth) In this example, we carry out a simulation of crystal growth with
fourfold anisotropy, and investigate 
how the anisotropic coefficient and the latent heat coefficient $K$ affect the shape of the dendritic crystal. 

We consider a benchmark problem, which has been extensively studied; see, e.g., \cite{KR99,K93modeling,ZCY19novel,Y21novel}. 
We set the initial conditions as follows:
\begin{align*}
\phi(x,y,0)=\tanh\Big(\frac{r_0-((x-x_0)^2+(y-y_0)^2)}{\varepsilon_0}\Big),
\quad   T(x,y,0)=\left\{
\begin{aligned}
0,\quad\   &  & \phi(x,y,0)>0,\\
-0.6, &  &  otherwise,
\end{aligned}
\right.
\end{align*}
where $r_0=9e-4,\ x_0=y_0=0,\ \varepsilon_0=1.8e-4$. The model parameters used in the simulation are:
\begin{align*}
&\varrho=1e3,\  \varepsilon=0.015,\  \sigma=0.1,\  \lambda=4e2,\    D=2.5e-3,\\
&S_1=0.6,\  S_2=10,\  S_3=S_4=4,\  B=4e5,\ \tau=0.01.
\end{align*}
We let the latent heat parameter $K$ vary. The spatial spectral discretization uses the polynomial space of degree $512$
at each direction. 

\begin{figure}[!ht]
\centering
\subfigure[$\phi$ at $t=0,\ 3,\ 6,\ 9$. $K=0.6$]{
	\begin{minipage}[t]{0.2\linewidth}
		\includegraphics[width=1.6in]{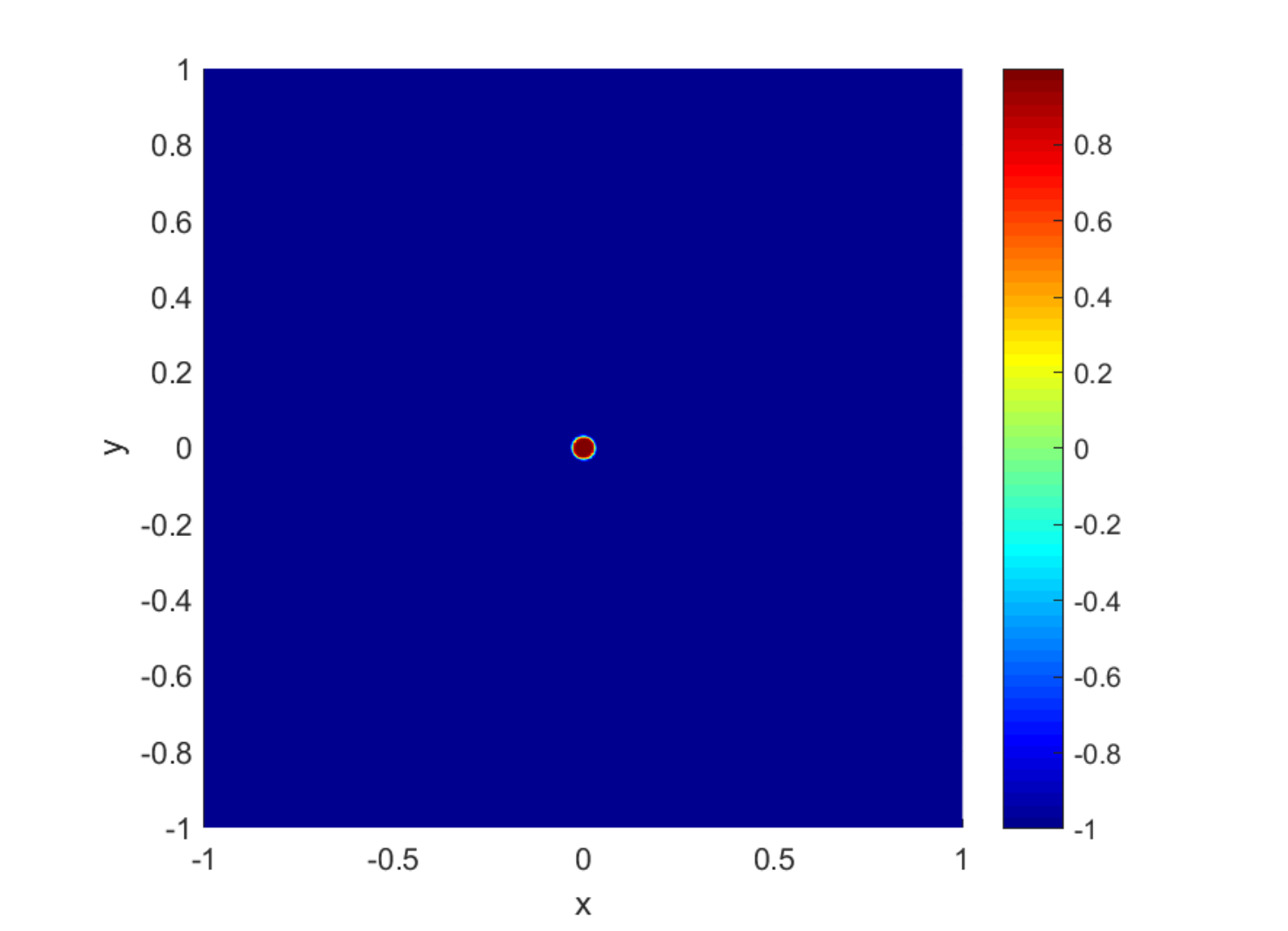}
	\end{minipage}
	\begin{minipage}[t]{0.2\linewidth}
		\includegraphics[width=1.6in]{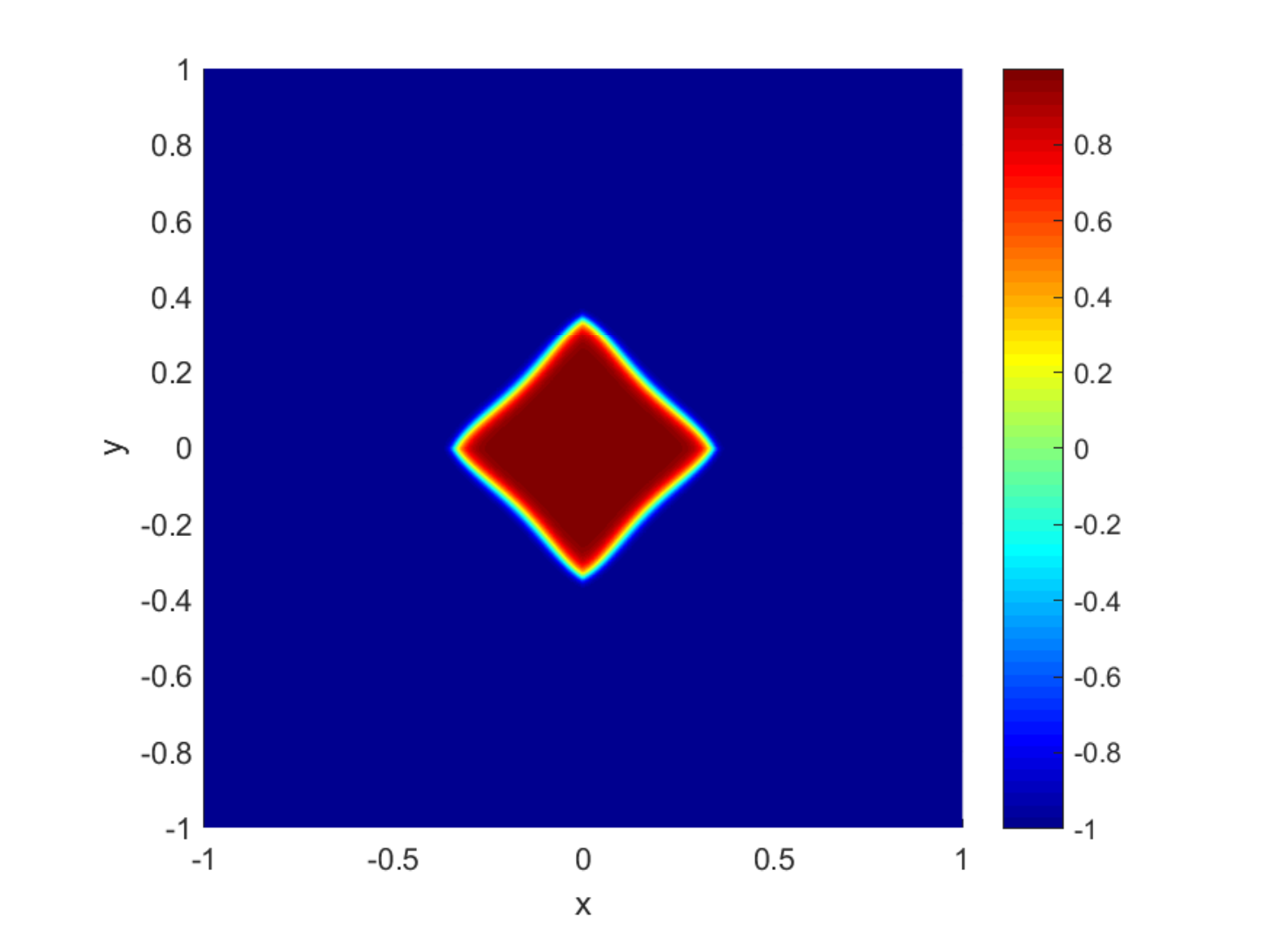}
	\end{minipage}
	\begin{minipage}[t]{0.2\linewidth}
		\includegraphics[width=1.6in]{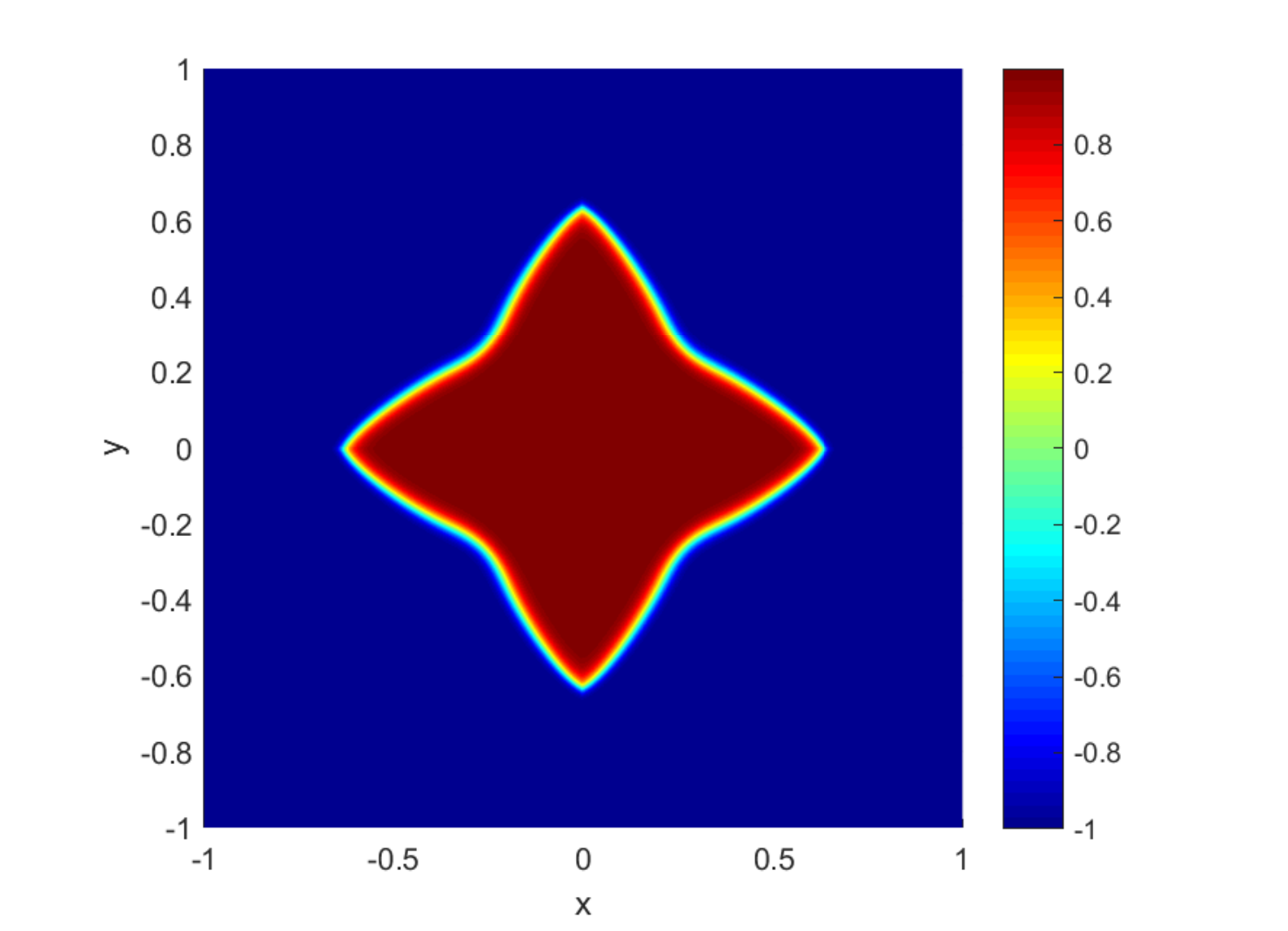}
	\end{minipage}
	\begin{minipage}[t]{0.2\linewidth}
		\includegraphics[width=1.6in]{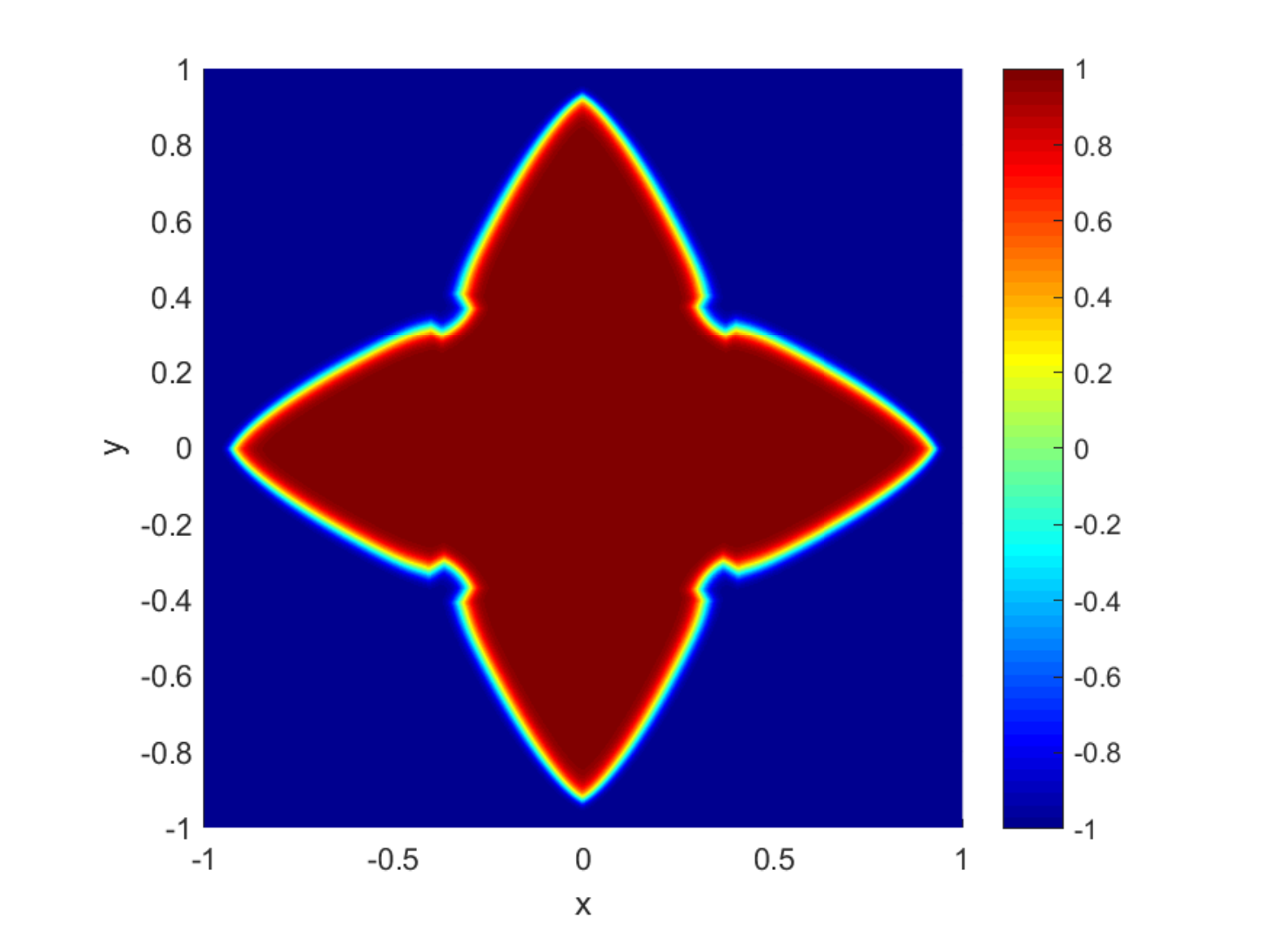}
	\end{minipage}
}
\quad
\subfigure[$\phi$ at $t=3,\ 6,\ 9,\ 11$. $K=0.8$]{
	\begin{minipage}[t]{0.2\linewidth}
		\includegraphics[width=1.6in]{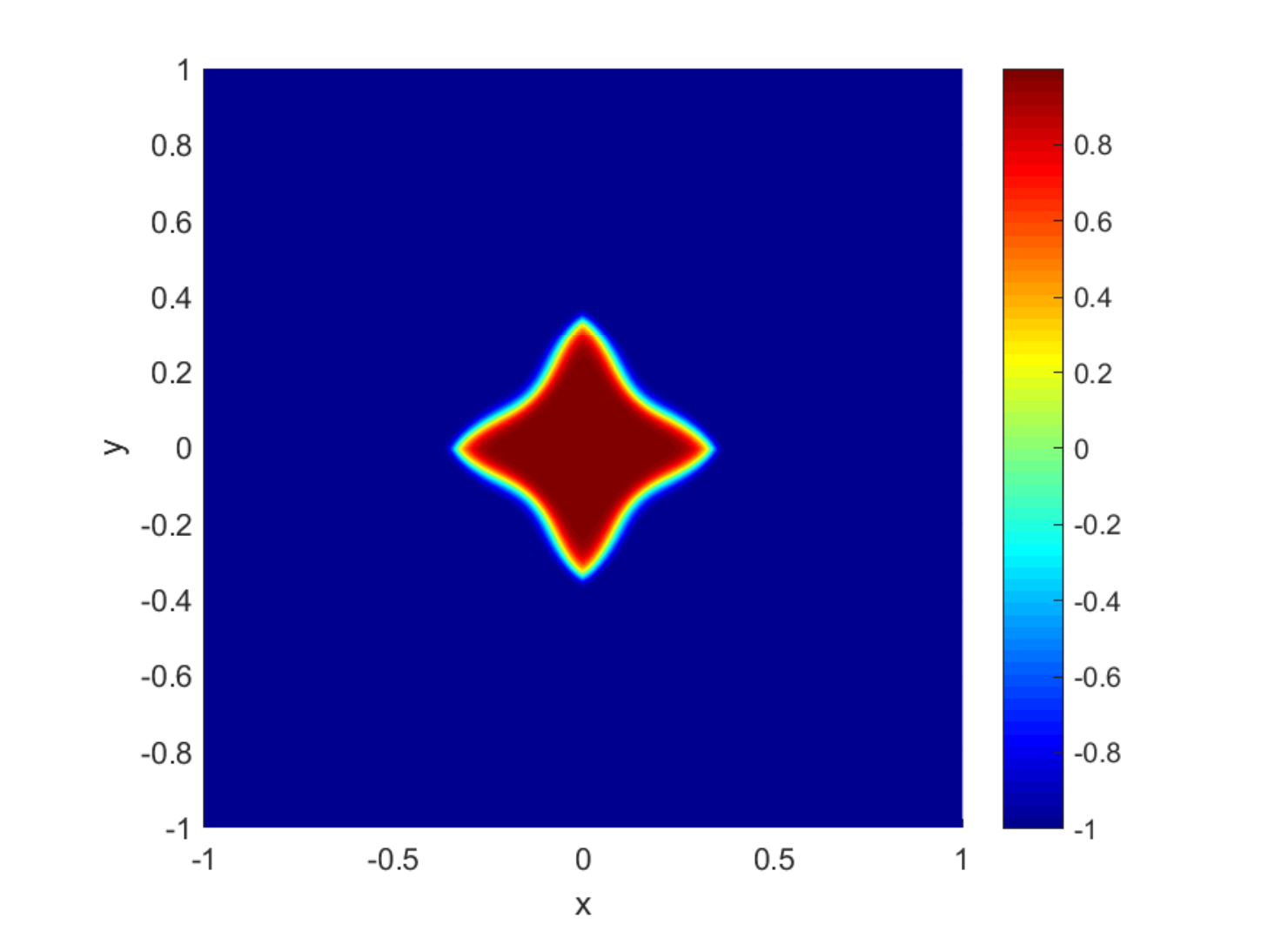}
	\end{minipage}
	\begin{minipage}[t]{0.2\linewidth}
		\includegraphics[width=1.6in]{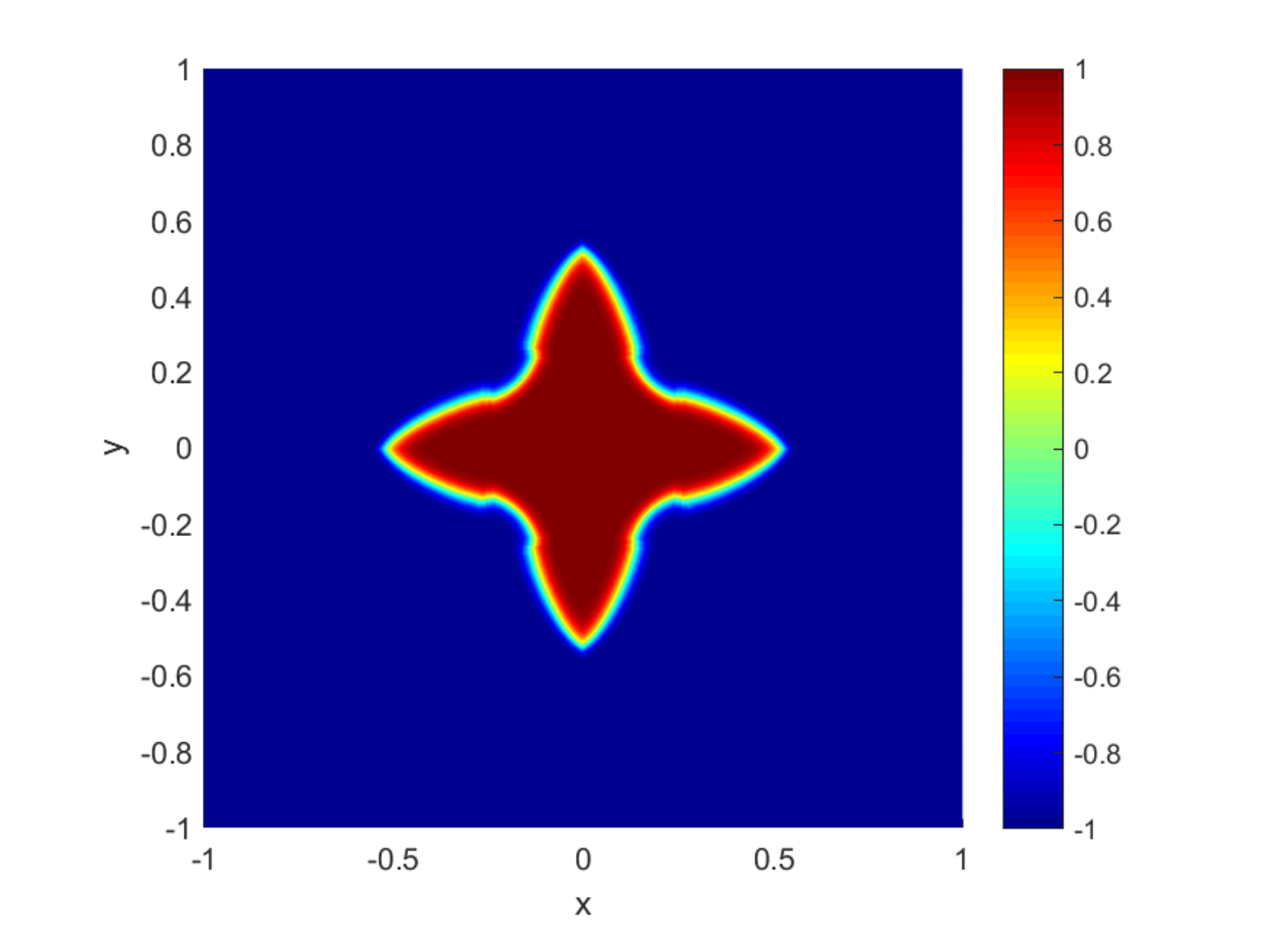}
	\end{minipage}
	\begin{minipage}[t]{0.2\linewidth}
		\includegraphics[width=1.6in]{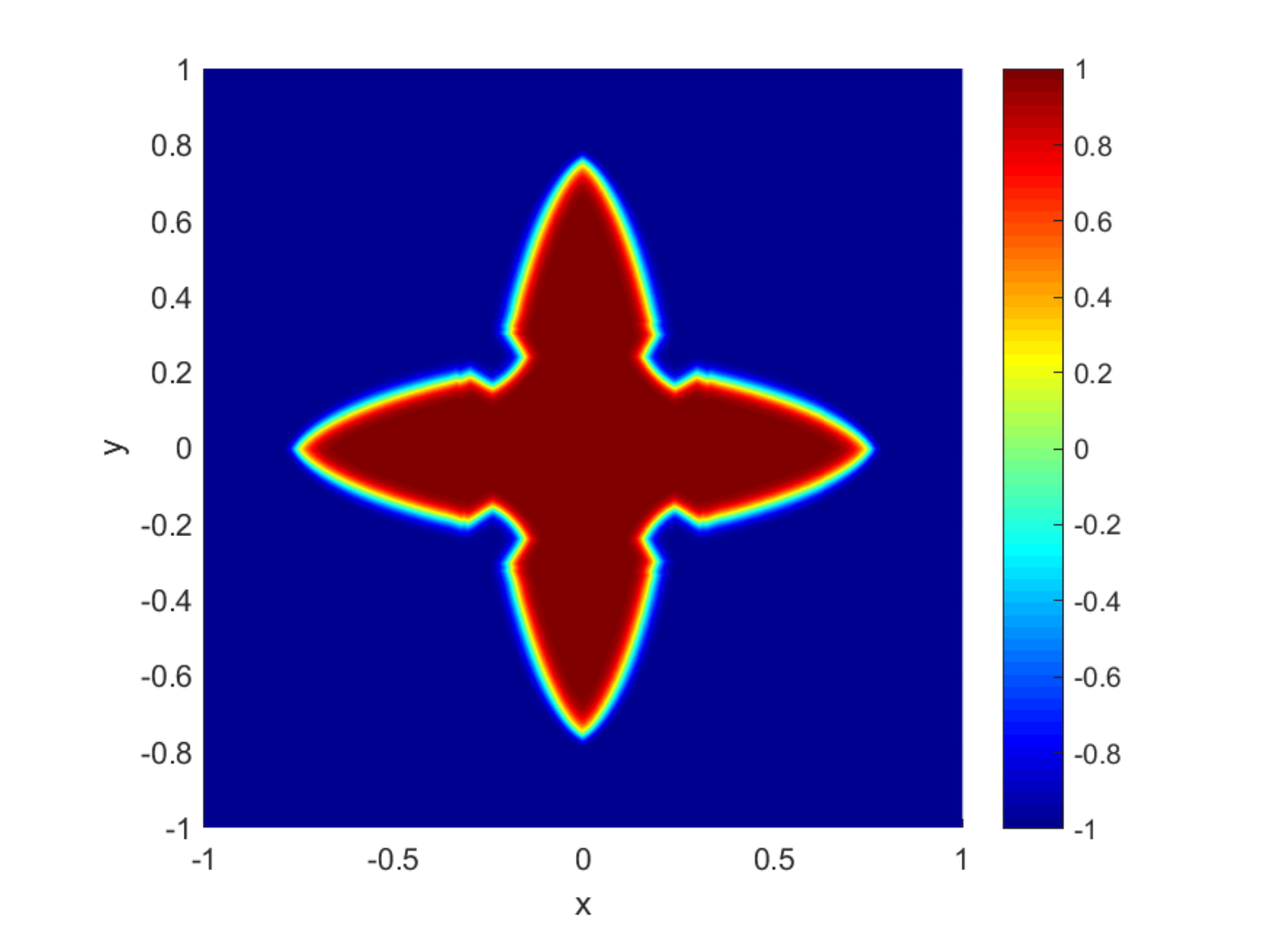}
	\end{minipage}
	\begin{minipage}[t]{0.2\linewidth}
		\includegraphics[width=1.6in]{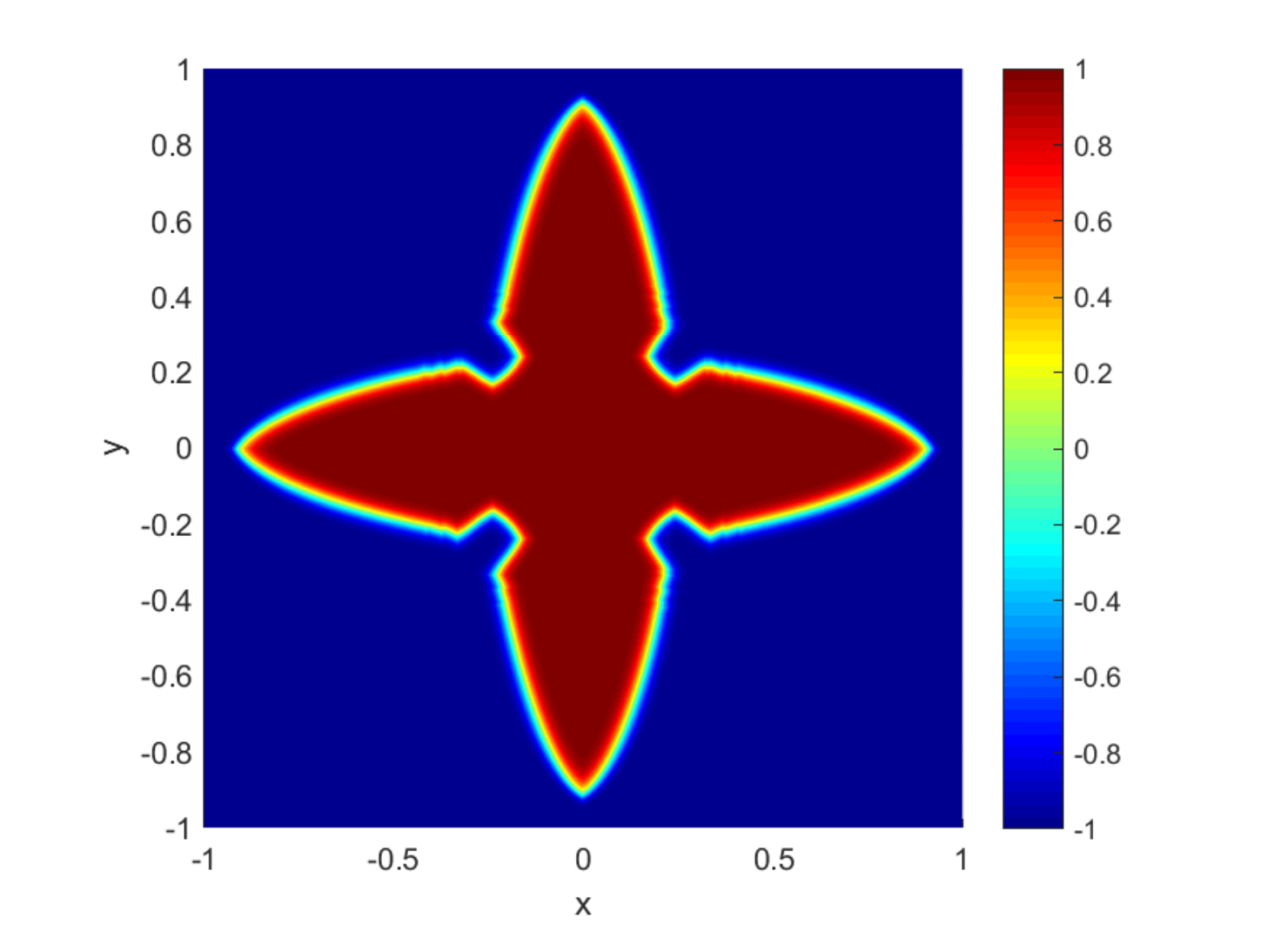}
	\end{minipage}
}%
\quad
\subfigure[$\phi$ at $t=6,\ 9,\ 11,\ 14$. $K=1$]{
	\begin{minipage}[t]{0.2\linewidth}
		\includegraphics[width=1.6in]{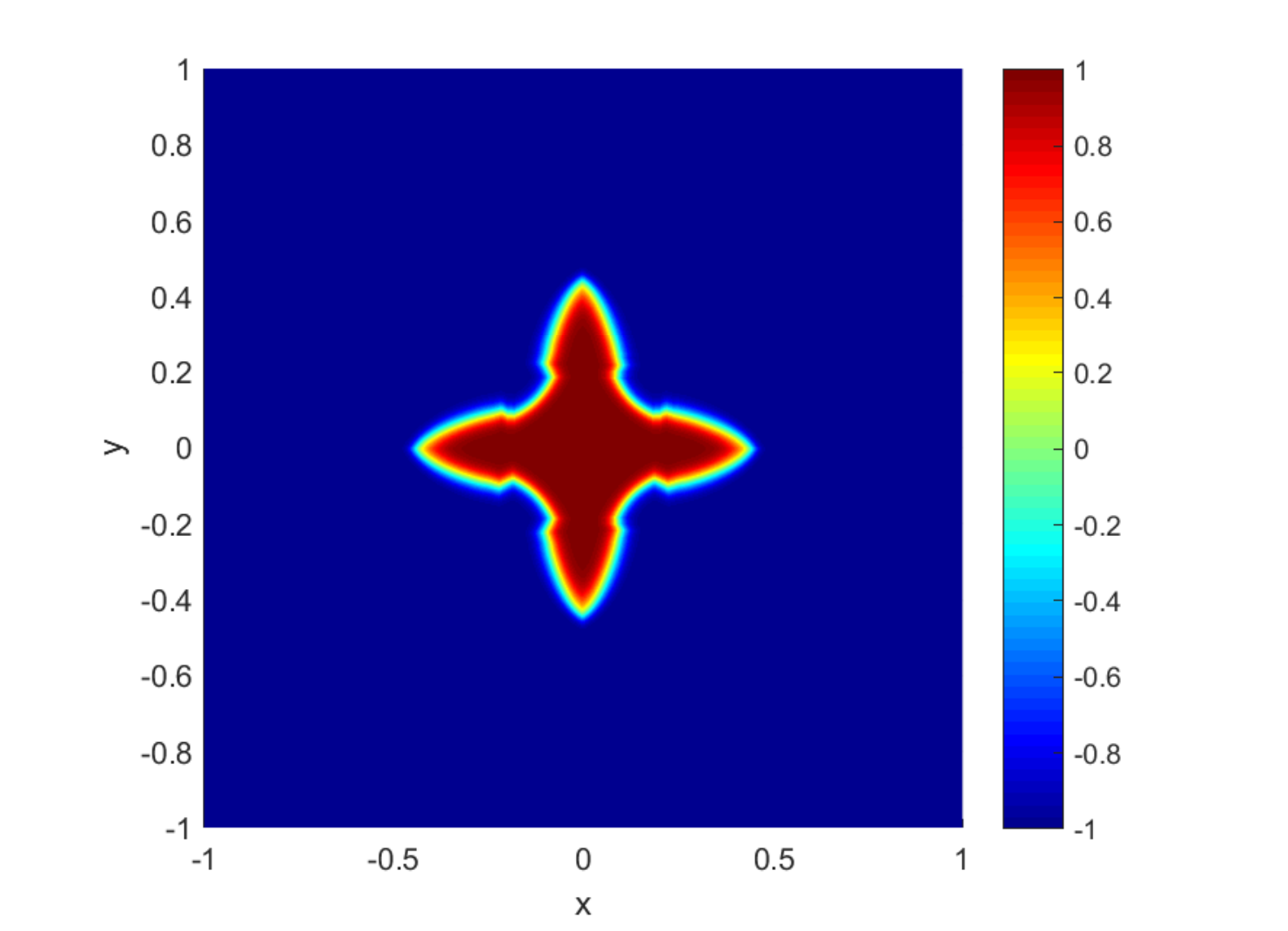}
	\end{minipage}
	\begin{minipage}[t]{0.2\linewidth}
		\includegraphics[width=1.6in]{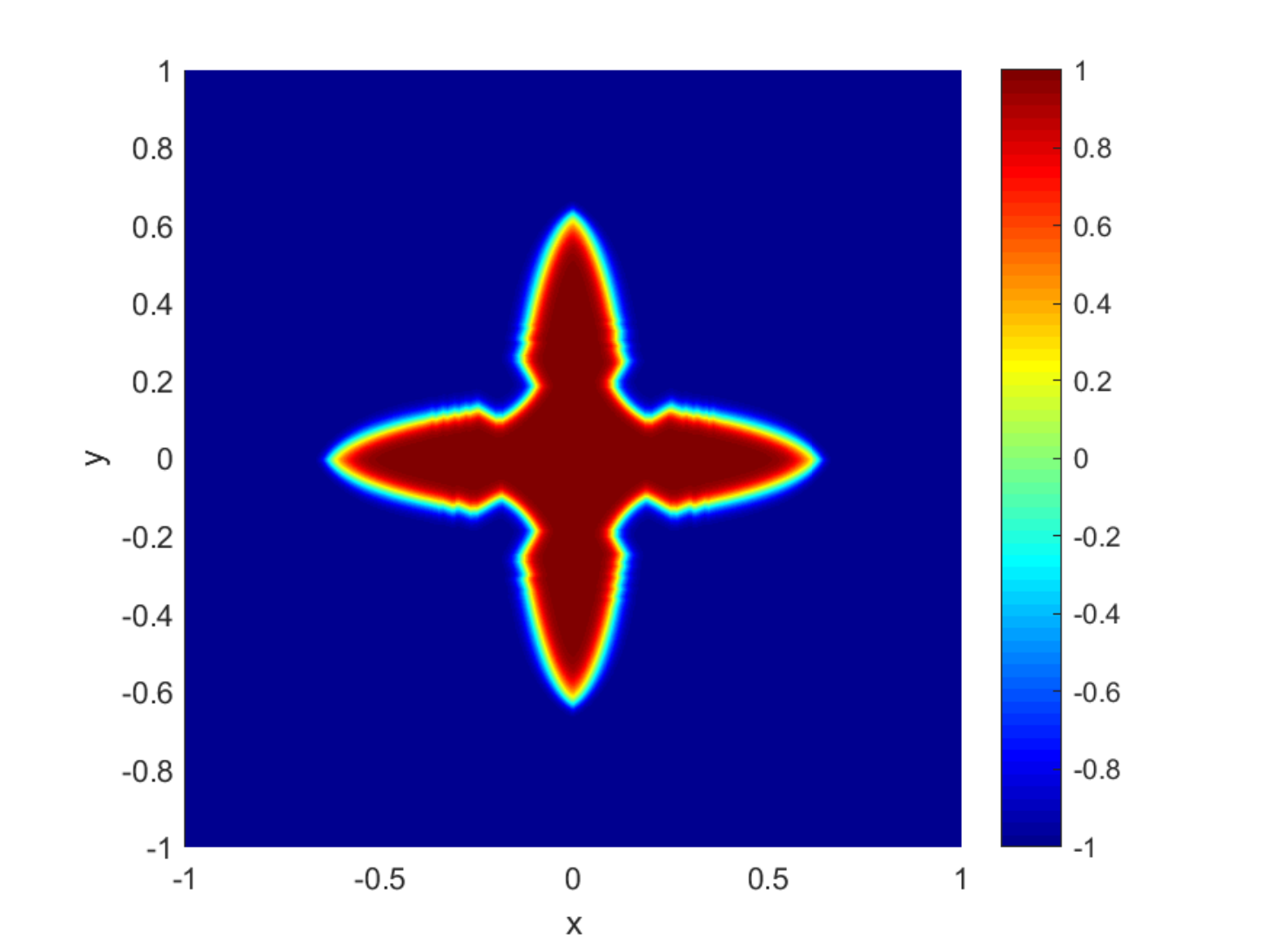}
	\end{minipage}
	\begin{minipage}[t]{0.2\linewidth}
		\includegraphics[width=1.6in]{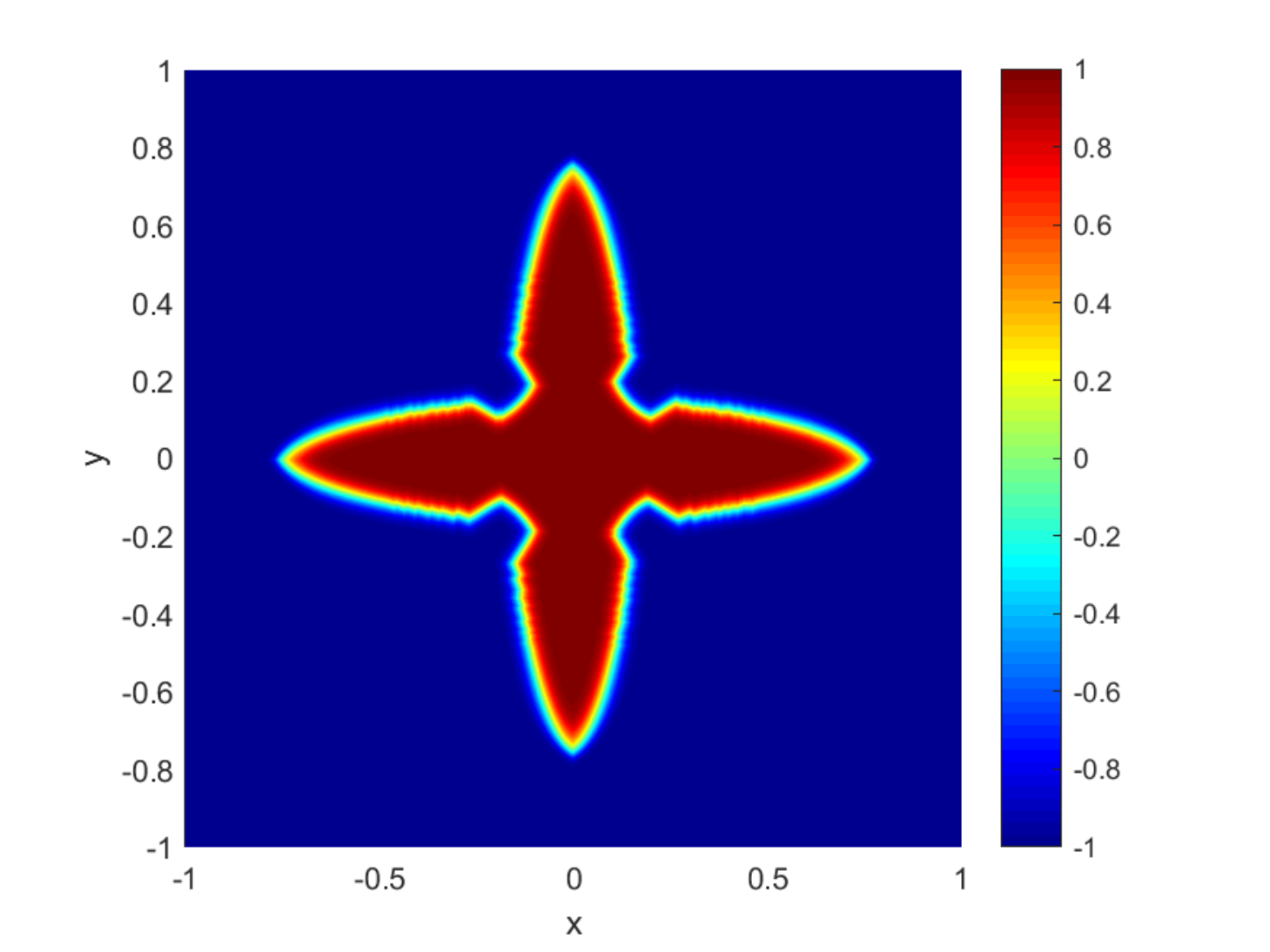}
	\end{minipage}
	\begin{minipage}[t]{0.2\linewidth}
		\includegraphics[width=1.6in]{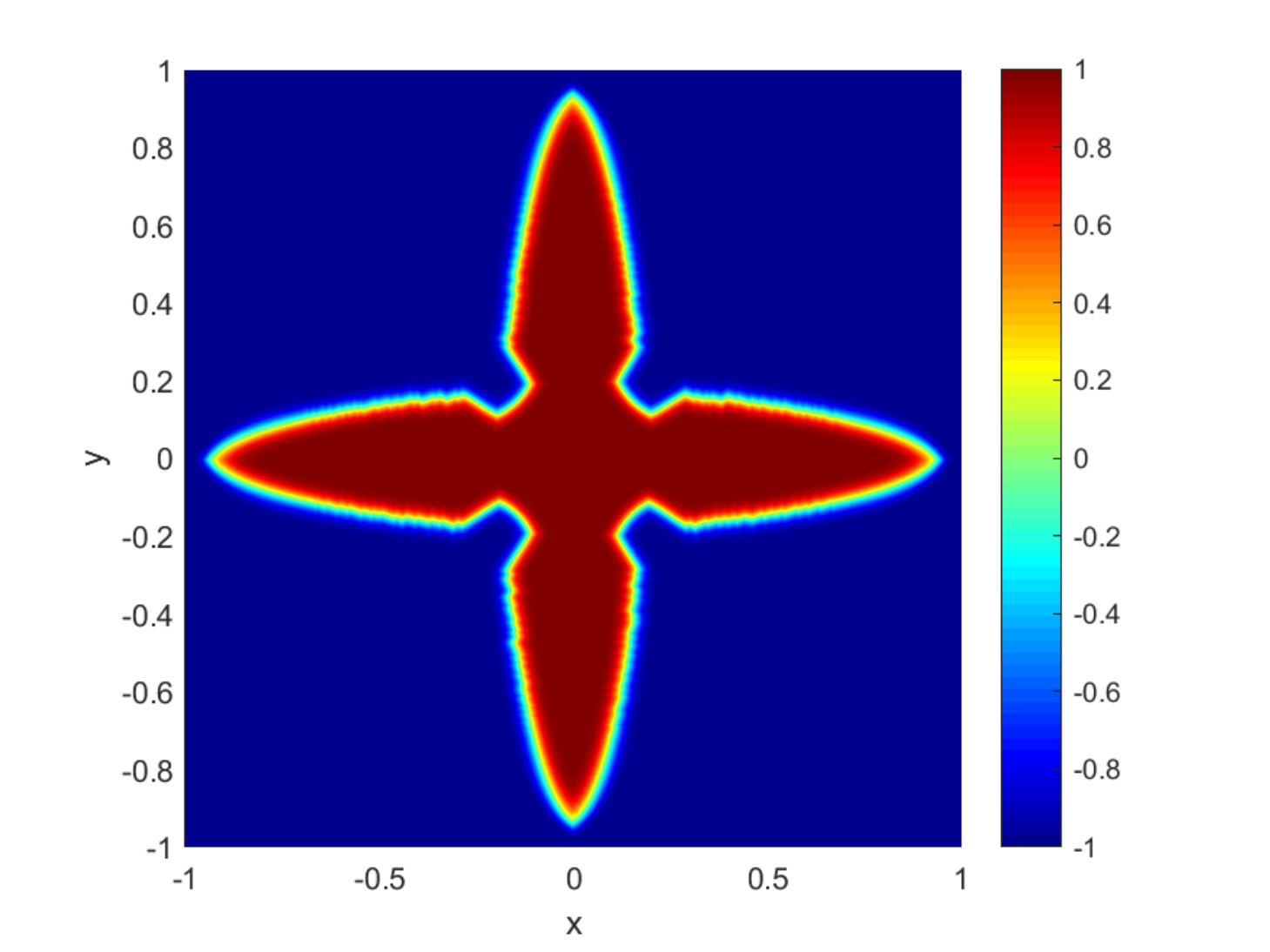}
	\end{minipage}
}%
\quad
\subfigure[$\phi$ at $t=9,\ 11,\ 14,\ 17$. $K=1.2$]{
	\begin{minipage}[t]{0.2\linewidth}
		\includegraphics[width=1.6in]{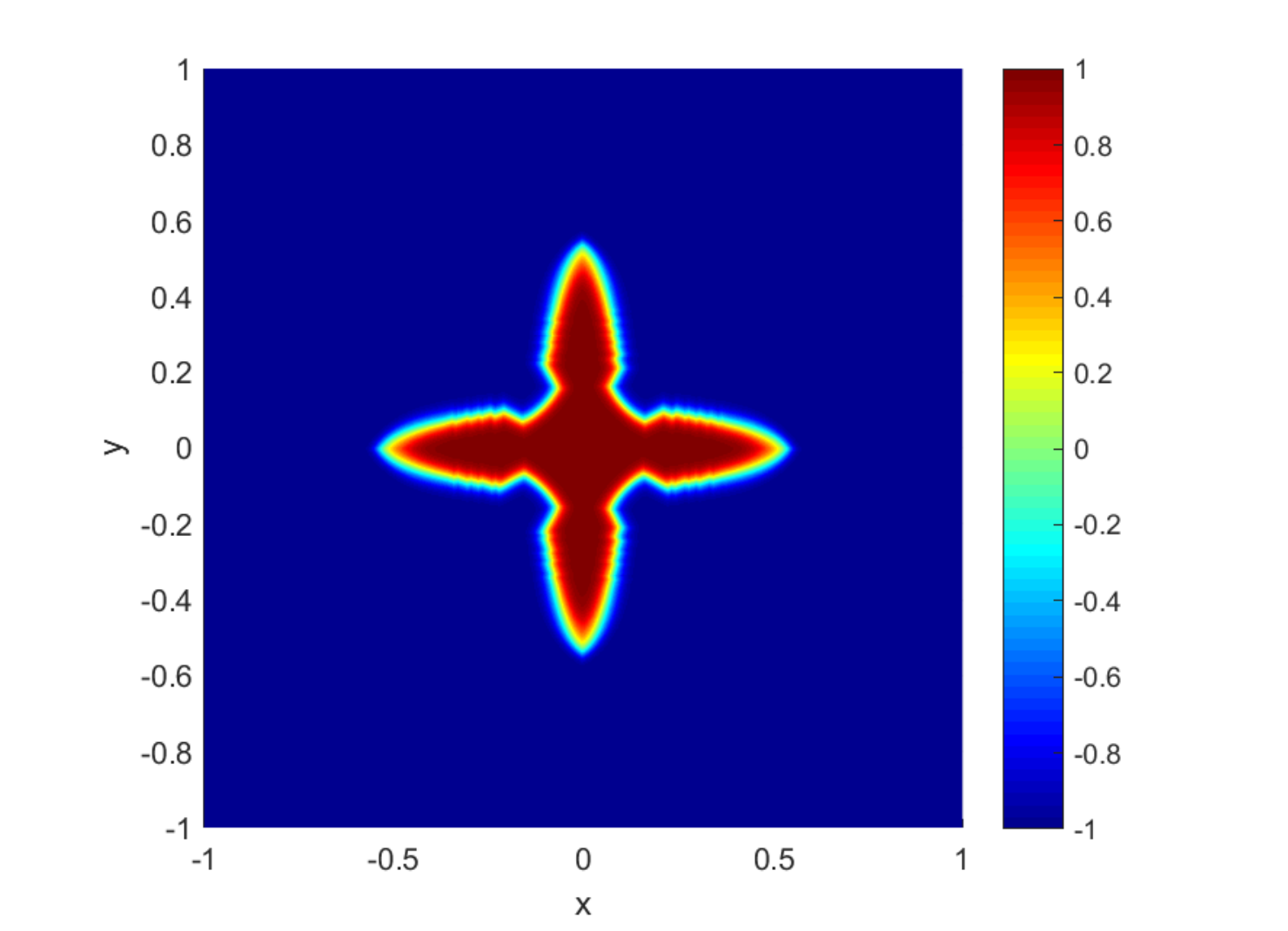}
	\end{minipage}
	\begin{minipage}[t]{0.2\linewidth}
		\includegraphics[width=1.6in]{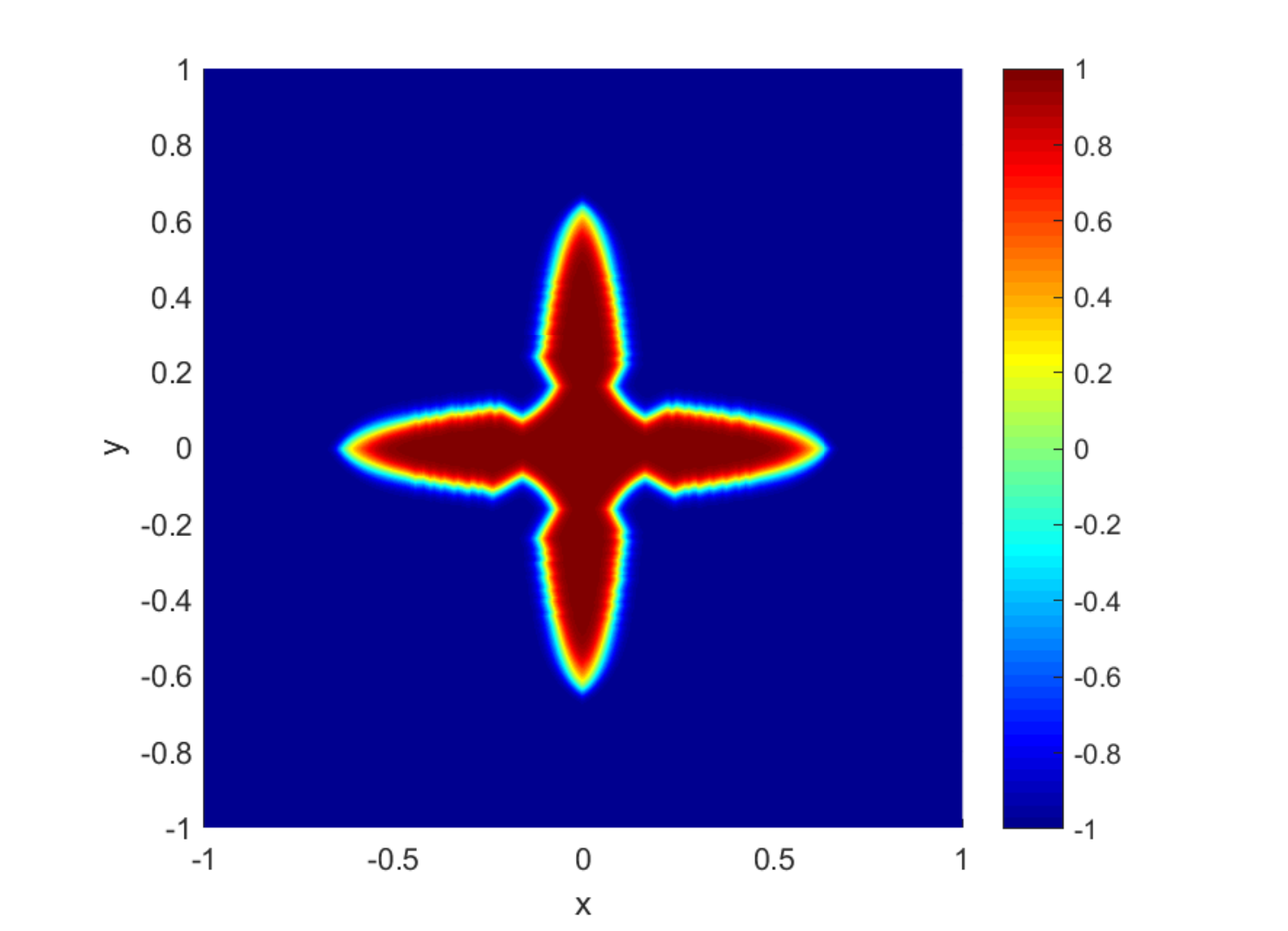}
	\end{minipage}
	\begin{minipage}[t]{0.2\linewidth}
		\includegraphics[width=1.6in]{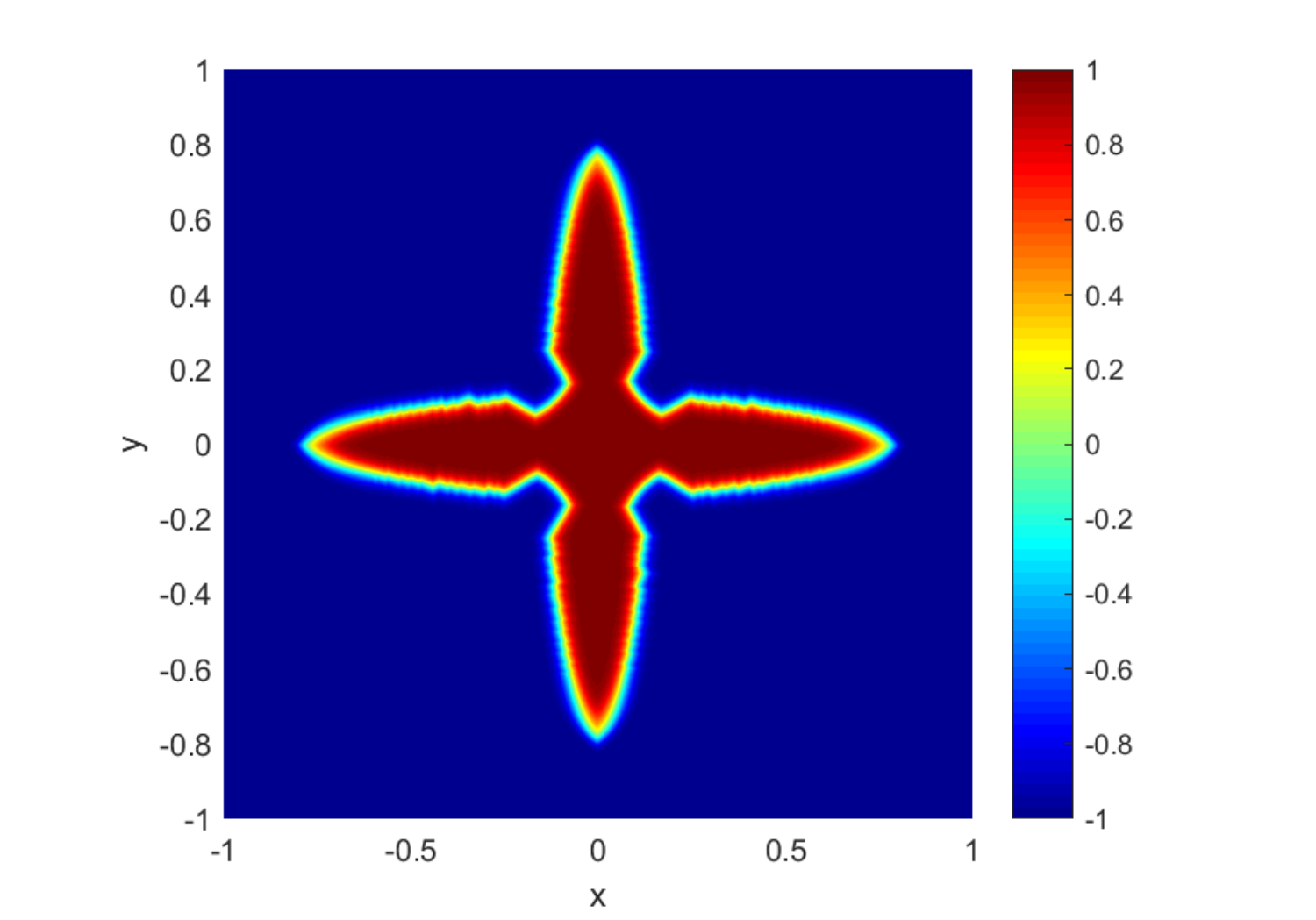}
	\end{minipage}
	\begin{minipage}[t]{0.2\linewidth}
		\includegraphics[width=1.6in]{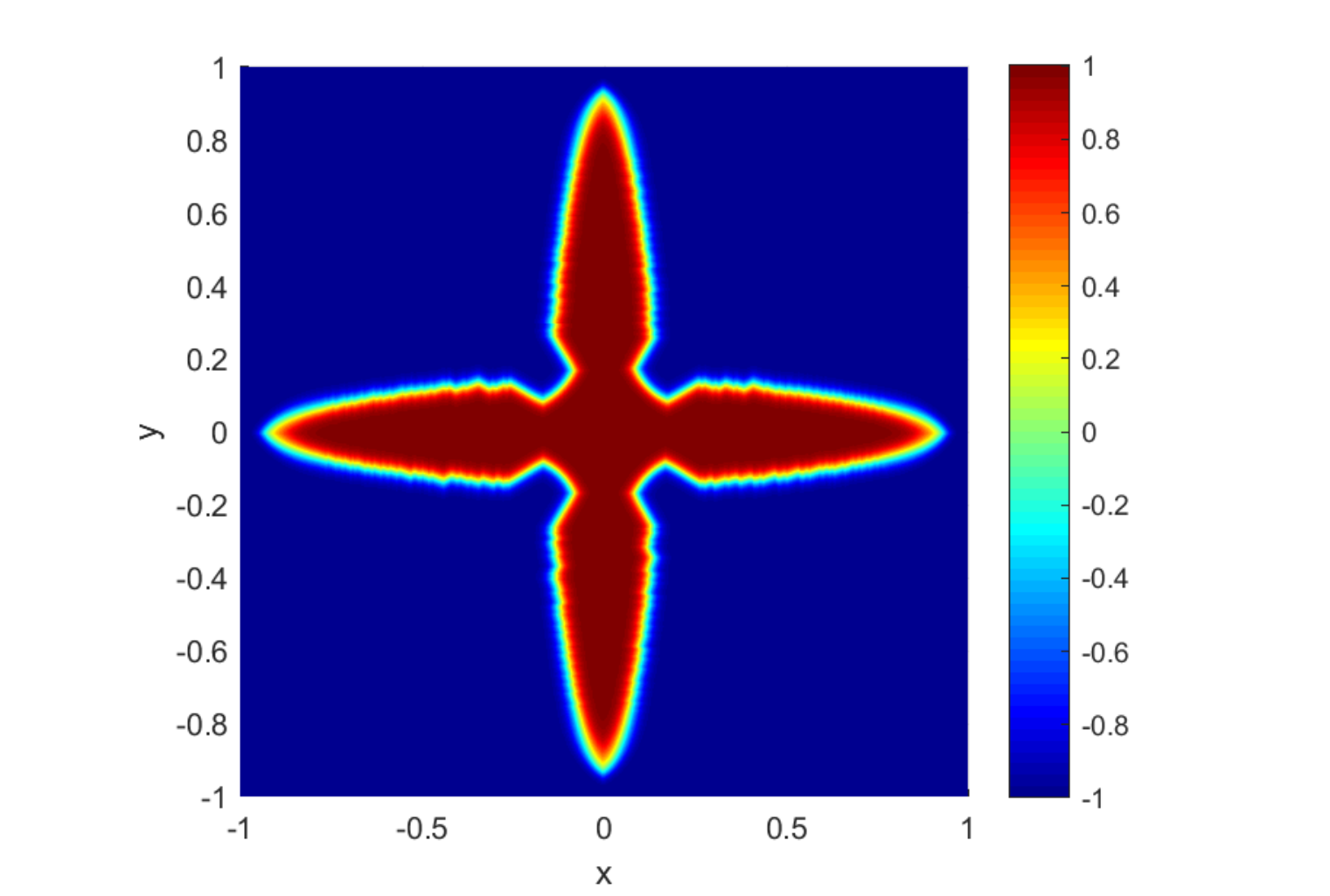}
	\end{minipage}
}
\quad
\subfigure[Temperature field $T$ at the last moment of above cases. 
From left to right: $K=0.6,\  K=0.8,\ K=1,\ K=1.2$.]{
	\begin{minipage}[t]{0.2\linewidth}
		\includegraphics[width=1.6in]{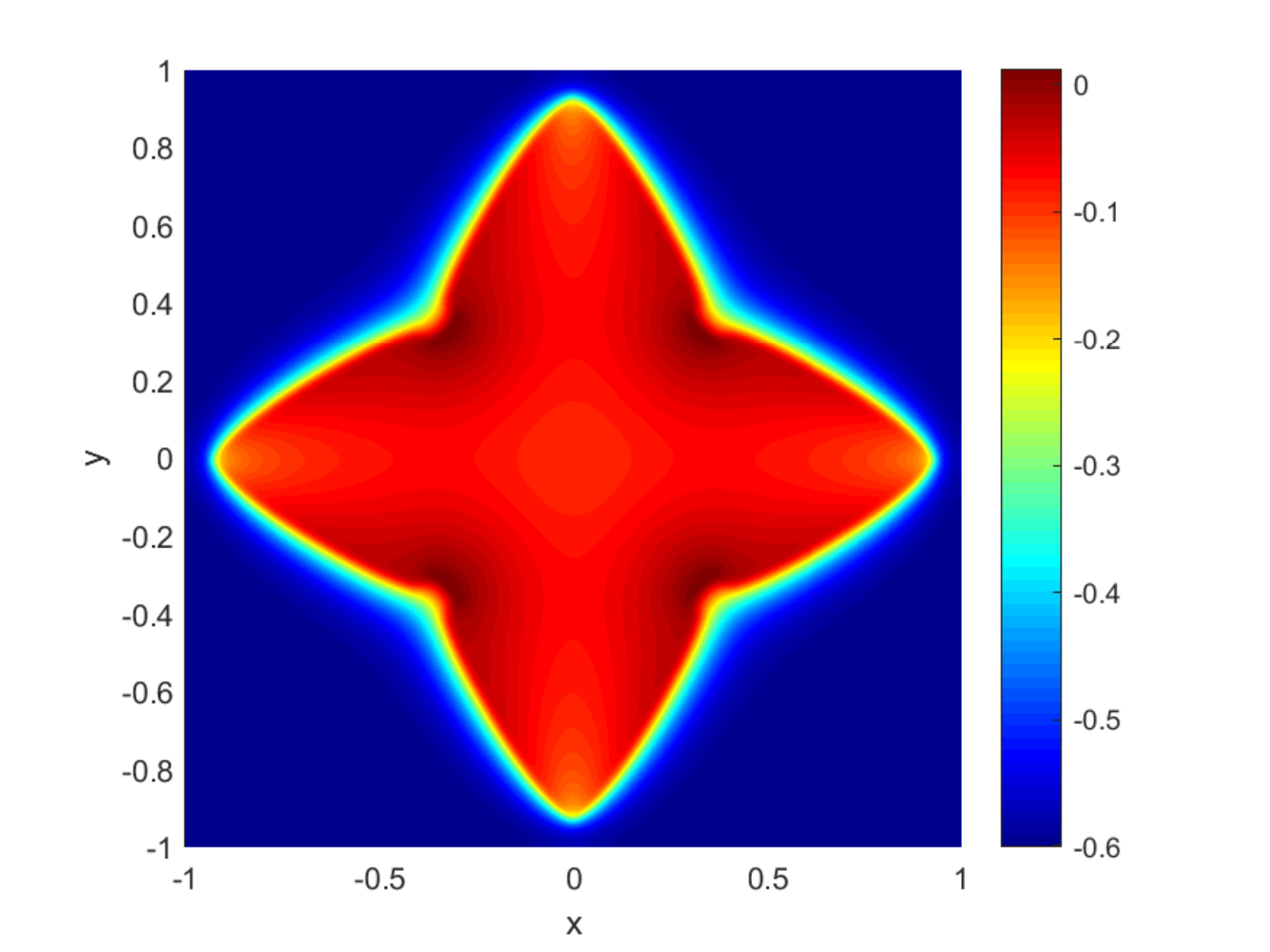}
	\end{minipage}
	\begin{minipage}[t]{0.2\linewidth}
		\includegraphics[width=1.6in]{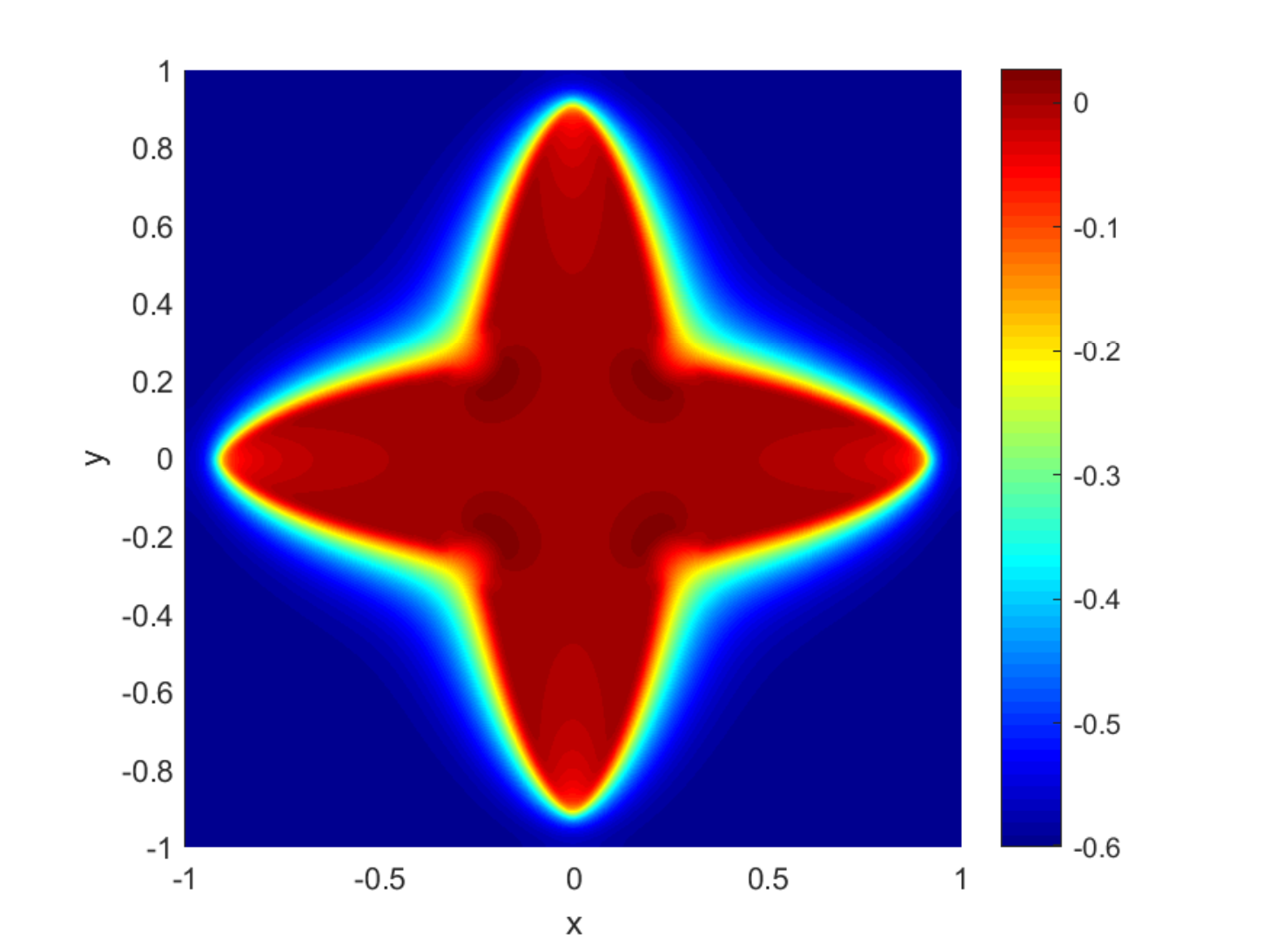}
	\end{minipage}
	\begin{minipage}[t]{0.2\linewidth}
		\includegraphics[width=1.6in]{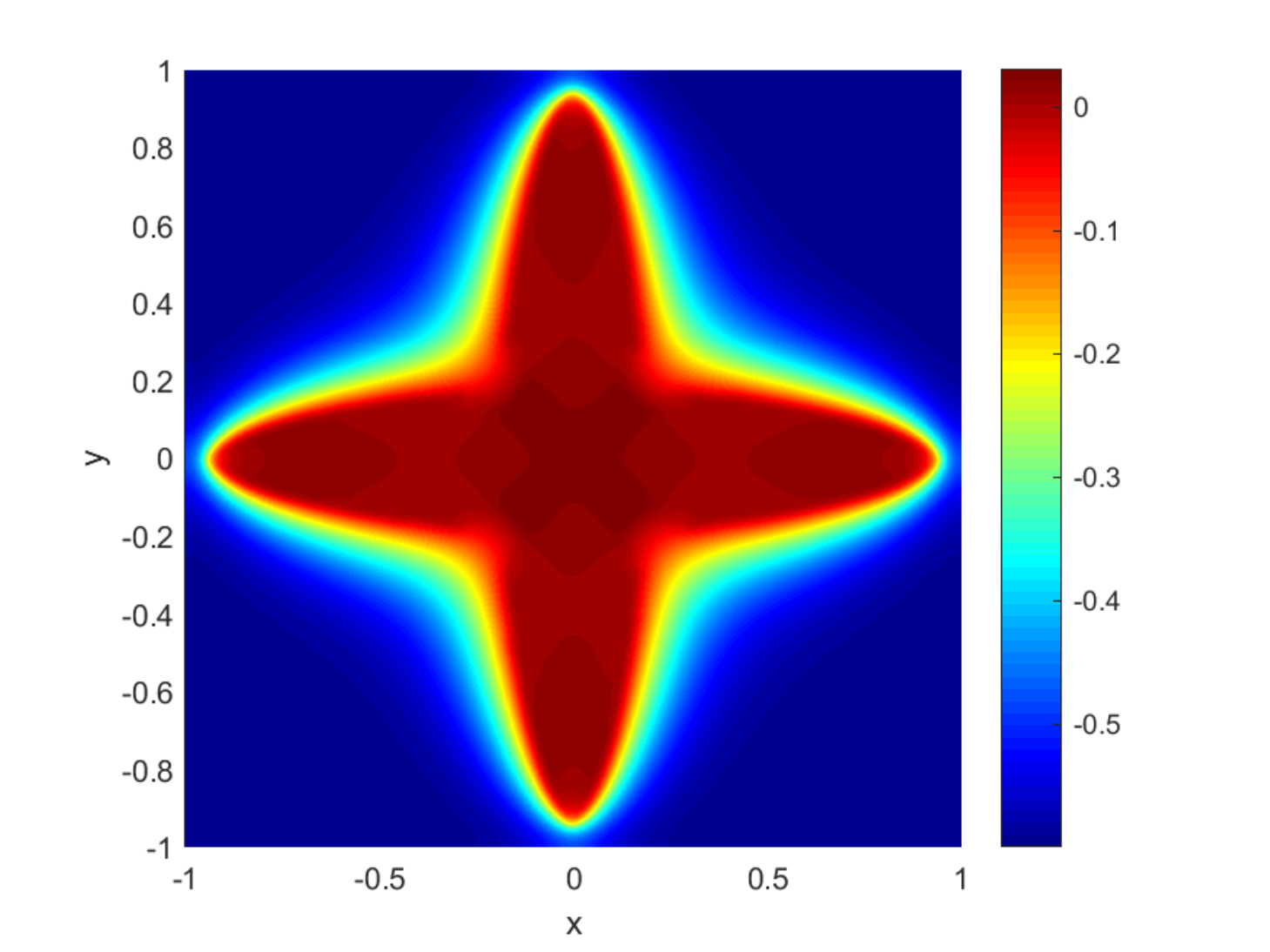}
	\end{minipage}
	\begin{minipage}[t]{0.2\linewidth}
		\includegraphics[width=1.6in]{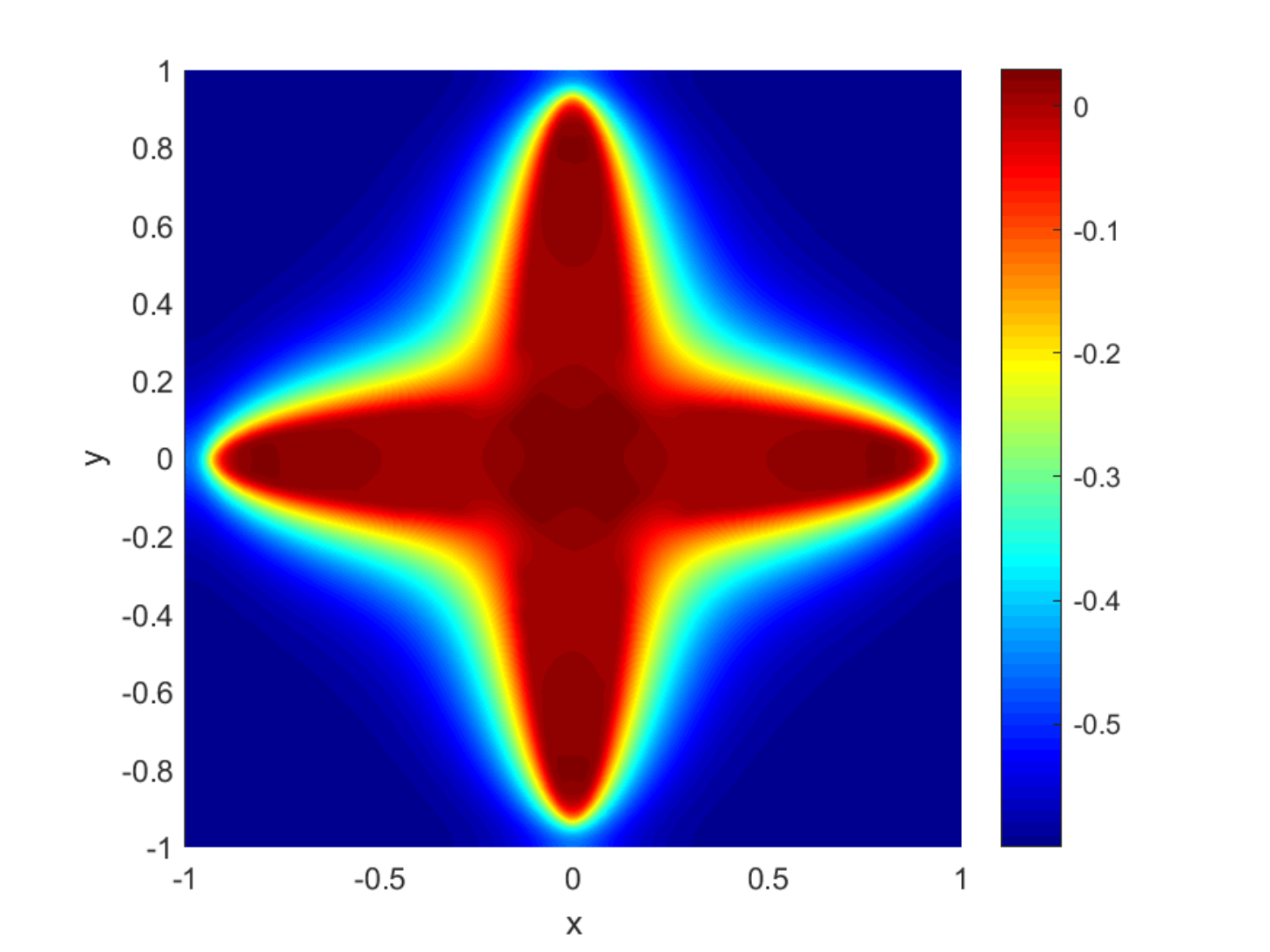}
	\end{minipage}
}
\centering
\caption{(Example \ref{exp3}) Dendritic crystal growth with fourfold anisotropy for different values of the latent heat parameter $K$. 
	(a)-(d): snapshots of the phase field $\phi$ at different times; (e): the temperature field $T$ at the last moment of (a)-(d). }
\label{figeg3}
\end{figure}

In Figures \ref{figeg3}(a)-(d), we present snapshots of $\phi$ at different time instances for $K$ varying from 0.6 to 1.2 
with an incremental value 0.2. The isocontours of the phase field function $\phi$ observed from the figures
clearly indicate that the four prominent branches are always formed in all cases starting with the same small circle.
Moreover, it is seen that the parameter $K$ affects the width of the branch: larger is $K$, thinner is the width of the branches, and 
sharper are the tips.
\begin{figure}[!ht]
\centering
\subfigure[$E^n$]{
	\includegraphics[width=7.3cm]{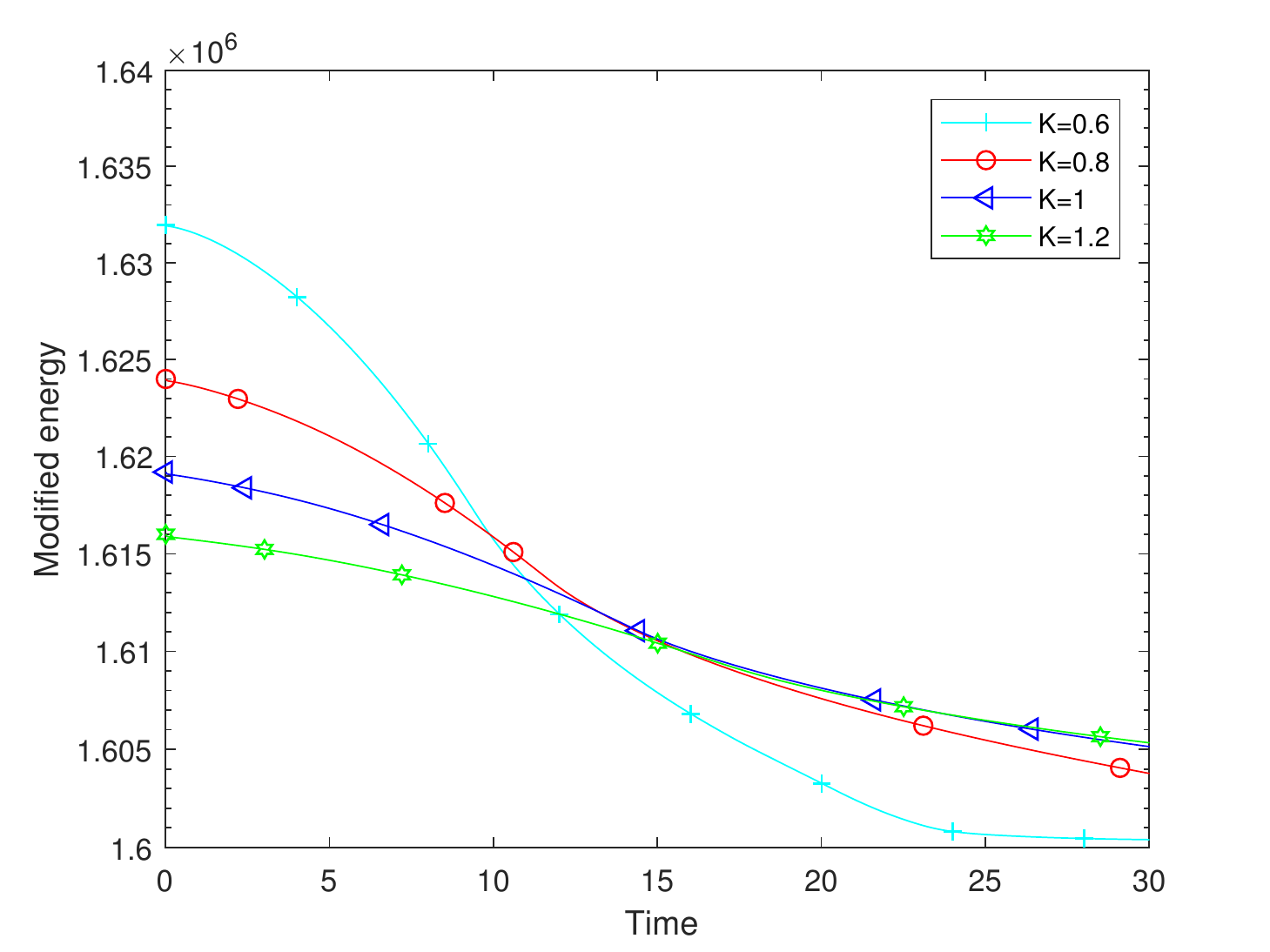}
}
\hspace{-10mm}
\subfigure[Crystal area $\int_{\Omega} \frac{1+\phi}{2}\ d\x$]{
	\includegraphics[width=7.3cm]{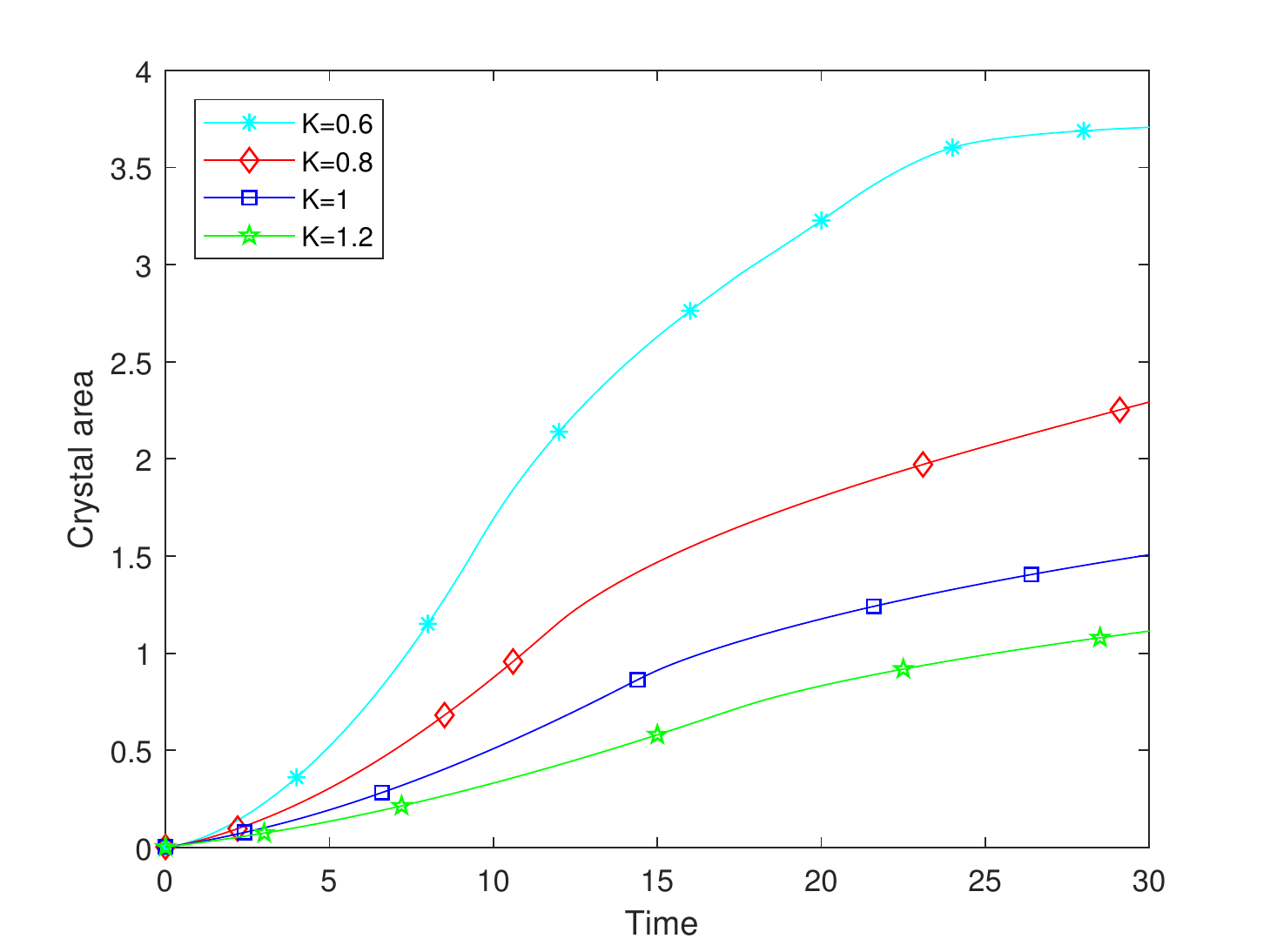}}
\caption{(Example \ref{exp3}) Time evolution of the modified energy $E^n$ and crystal area 
	$\int_{\Omega} \frac{1+\phi}{2}\ d\x$ for different latent heat parameter $K$.}
\label{figEnA}
\end{figure}

The isocontours of the temperature $T$ at the last moment of each simulation is plotted in
the figures \ref{figeg3}(e). It is observed that the contours of the temperature $T$ take similar dendrite crystal shape as
the phase field. This is due to the fact that the heat is propagating only at the interface.

The dissipation behavior of the modified energy $E^n$ in time is shown in Figure \ref{figEnA}(a).
The monotonic decay feature of $E^n$ for all tested $K$ reflects good stability property of the scheme used in the calculation.  
Finally we give in Figure \ref{figEnA}(b) evolution of the area of the crystal, defined by the quantity $\int_{\Omega}\frac{1+\phi}{2} d\x$,
for several values of $K$.
We see that the area of the crystal keeps increasing during the simulation. 
This is in a good agreement with the existing results; see, e.g., 
\cite{KR99, Y19efficient, Y21anovel}. 

\end{example}

\section{Concluding  remarks}
We have proposed a class of new time-stepping schemes for the anisotropic phase-field dendritic crystal growth model. 
The proposed schemes were constructed 
based on an auxiliary variable approach for the Allen-Cahn equation and sophisticated treatment
of the terms coupling the Allen-Cahn equation and temperature equation. 
In particular, the new reformulation of the model introduced 
in the paper plays a key role in developing efficient schemes.
Thanks to the carefully chosen extra terms added to 
the time discretization,
we were able to construct a second-order scheme, which is linear, 
decoupled, uniquely solvable, and unconditionally stable. 
A detailed comparison with existing schemes is given, and the
advantage of the new schemes are emphasized. 
The stability property of the proposed schemes was rigorously established, while the  
convergence rate was carefully examined through
a series of numerical tests.
Our analysis and numerical experiments demonstrated 
the efficiency of the proposed method.
It seems to us that the approach proposed in this paper 
is extendable to more complex models such as those studied in \cite{Y21new,Y21novel,Y21numerical}. 

\bibliographystyle{siam}
\bibliography{mybibfile}
\end{document}